\documentclass[a4paper,12pt]{amsart}
\usepackage{amssymb}
\usepackage{amscd}
\usepackage{extarrows}
\usepackage[colorlinks=true]{hyperref}
\usepackage{tikz}
\usepackage{enumerate}
\usetikzlibrary{arrows,matrix}
\usetikzlibrary{shapes.geometric}
\allowdisplaybreaks[1]%allows to break equations if really necessary
\long\def\rem#1{}
\title{Path representation of maximal parabolic Kazhdan--Lusztig polynomials}
\author[K.~Shigechi]{Keiichi~Shigechi}
\address{K.~Shigechi, Institute for Theoretical Physics, Valckenierstraat 65, 1018 XE Amsterdam, The Netherlands}
\email{k.shigechi@uva.nl}
\thanks{KS's work is part of the research programme of the Foundation for Fundamental Research on Matter (FOM), which is financially supported by the Netherlands Organisation for Scientific Research (NWO).}
\author[P.~Zinn-Justin]{Paul~Zinn-Justin}
\address{P.~Zinn-Justin, UPMC Univ Paris 6, CNRS UMR 7589, LPTHE,
75252 Paris Cedex, France}
\email{pzinn@lpthe.jussieu.fr}
\thanks{PZJ acknowledges partial support
from ESF program ``MISGAM'',
ANR program ``GRANMA'' BLAN08-1-13695,
and KNAW's visiting professors programme.}
\date{\today}
\newcommand\R{R}

\newcommand\Z{\mathbb{Z}}
\renewcommand\S{\mathcal{S}}
\renewcommand\H{\mathcal{H}_N}
\newcommand\M{\mathcal{M}_{N,K}}
\newcommand\Md{\mathcal{M}_{N,N-K}}
\newcommand\coset{\S_N/(\S_K\times\S_{N-K})}
\newcommand\cosetb{\S_N/(\S_{N-K}\times\S_{K})}
\newcommand\short[1]{{\rm s}(#1)}
\newcommand\Card{\#}
\newtheorem{example}{Example}
\newtheorem{lemma}{Lemma}
\newtheorem{theorem}{Theorem}
\newtheorem{defn}{Definition}
\newtheorem{prop}{Proposition}
\newtheorem{cor}{Corollary}

\newtheorem{remark}{Remark}
\makeatletter
\newcommand{\gettikzxy}[3]{% returns coordinates expressed *as dimensions*
  \tikz@scan@one@point\pgfutil@firstofone#1\relax
  \edef#2{\the\pgf@x}%
  \edef#3{\the\pgf@y}%
}
\makeatother
\newdimen\linkpatternunit%
\newcount\linkpatternsize%
\newcount\lpsize%the global internal variable... ugly but \foreach forces us to use a global variable
\newcount\linkpatternfused%for fused link patterns, level of fusion
\newif\iflinkpatterninverted%mirror symmetry
\newif\iflinkpatternhalfnumbered%only used by \tanglelinkpattern = only lower #
\newif\iflinkpatternnumbered%labels = numbers
\newif\iflinkpatternrenumbered%labels = redefined arbitrarily
\newif\iflinkpatterntikzstarted%
\newif\iflinkpatternboxed%draw a box around link patterns/tangles
\newif\iflinkpatternaxis%draw a line at starting/endpoints
\pgfkeys{/linkpattern/.cd,inverted/.is if=linkpatterninverted,numbered/.is if=linkpatternnumbered,tikzstarted/.is if=linkpatterntikzstarted,halfnumbered/.is if=linkpatternhalfnumbered,boxed/.is if=linkpatternboxed,axis/.is if=linkpatternaxis,vertexcolor/.store in=\linkpatternvertexcolor,edgecolor/.store in=\linkpatternedgecolor,tikzoptions/.style={every linkpattern/.append style={#1}},size/.code={\linkpatternsize=#1},numbering/.code={\linkpatternnumberedtrue\linkpatternrenumberedtrue\def\linkpatternnumbering{#1}},fused/.code={\linkpatternfused=#1},unit/.code={\linkpatternunit=#1},unpaired length/.store in=\linkpatternunpairedlength,height/.store in=\linkpatterntangleheight,every linkpattern/.style={x=\linkpatternunit,y=\linkpatternunit,baseline=0}}%change dir to /linkpattern, then list of optional arguments
\linkpatterninvertedfalse%
\linkpatternnumberedfalse%
\linkpatternrenumberedfalse%
\linkpatternhalfnumberedfalse%
\linkpatterntikzstartedfalse%
\linkpatternboxedfalse%
\linkpatternaxistrue%
\def\linkpatternvertexcolor{red}%
\def\linkpatternedgecolor{blue}%
\def\linkpatternunpairedlength{0.5}%unpaired min size
\def\linkpatterntangleheight{\optimaltangleheight}%height of tangle
\tikzset{vertex/.style={circle,thin,draw=black,fill=\linkpatternvertexcolor,inner sep=1.5pt}}%
\tikzset{edge/.style={very thick,draw=\linkpatternedgecolor}}%
\newcommand\linkpattern[2][]{%
{%we want the macro modifications to be local only
\pgfkeys{/linkpattern/.cd,#1}%parse the list of optional arguments
\iflinkpatternnumbered%
\tikzset{lbl/.style={label=\iflinkpatterninverted above\else below\fi:{$\scriptstyle ##1$}}}%
\else%
\tikzset{lbl/.style={}}
\fi%
\iflinkpatterntikzstarted\else%
\begin{tikzpicture}[/linkpattern/every linkpattern]%
\fi%
\iflinkpatterninverted%
\begin{scope}[yscale=-1]%
\fi%
\global\lpsize=\linkpatternsize
\foreach \x/\y in {#2}
{
\ifnum\x>\lpsize\global\lpsize=\x\fi
\ifnum\y>\lpsize\global\lpsize=\y\fi
}
\iflinkpatternboxed
\draw (0.5,0) rectangle ($(\lpsize+0.5,1)$);
\else
\iflinkpatternaxis
\draw (1,0) -- (\lpsize,0);
\else
\fi
\fi
\iflinkpatternrenumbered
\else
\def\linkpatternnumbering{1,...,\lpsize}
\fi
\foreach\x/\xx in \linkpatternnumbering
\coordinate[vertex,lbl={\xx}] (v\x) at (\x,0);
\foreach \x/\y/\z in {#2}
{
\ifx\y\z
\draw[bend left=60,edge]
\else
\draw[bend left=60,decoration={markings,mark = at position 0.5 with { \arrow[semithick]{\z} }},postaction={decorate},edge]
\fi
\ifnum\x<\y
(v\x) to (v\y);
\else
\ifnum\x>\y
(v\y) to (v\x);
\else
;%(v\x) -- ++(0,\linkpatternunpairedlength);
\fi
\fi
}
\iflinkpatterninverted
\end{scope}
\fi
\iflinkpatterntikzstarted\else%
\end{tikzpicture}%
\fi%
}}%
\newcount\maxsep
\def\linkpatternlooseness{0.6}% affects tangles only. may need to be adjusted....
\newcommand\tangle[2][]{%
{%we want the macro modifications to be local only
\pgfkeys{/linkpattern/.cd,#1}
\def\primetest##1'{}%
\def\hasaprime##1{\expandafter\primetest##1''}%aha!
\def\internalwithoutprime##1'{##1}%
\def\withoutprime##1{\if\hasaprime##1 %
\expandafter\internalwithoutprime##1\else ##1\fi}%
\iflinkpatternnumbered%
\iflinkpatternhalfnumbered%
\tikzset{lbla/.style={}}
\else%
\tikzset{lbla/.style={label=\iflinkpatterninverted below\else above\fi:{$\scriptstyle ##1$}}}%
\fi
\tikzset{lblb/.style={label=\iflinkpatterninverted above\else below\fi:{$\scriptstyle ##1$}}}%
\else%
\tikzset{lbla/.style={}}%
\tikzset{lblb/.style={}}%
\fi%
\iflinkpatterntikzstarted\else%
\begin{tikzpicture}[/linkpattern/every linkpattern]%
\fi%
\iflinkpatterninverted%
\begin{scope}[yscale=-1]%
\fi%
\global\lpsize=\linkpatternsize%
\newcount\tempx\newcount\tempxx%
\global\maxsep=0%
\foreach \x/\xx in {#2}%
{%
\tempx=\withoutprime{\x}
\tempxx=\withoutprime{\xx}
\ifnum\lpsize<\tempx\global\lpsize=\tempx\fi
\ifnum\lpsize<\tempxx\global\lpsize=\tempxx\fi%
\pgfmathsetcount{\maxsep}{max(\maxsep,abs(\tempx-\tempxx))}
\global\maxsep=\maxsep
}%
\pgfmathsetmacro{\optimaltangleheight}{0.6+0.4*\maxsep}
\iflinkpatternboxed
\draw (0.5,-0.5*\linkpatterntangleheight) rectangle (\lpsize+0.5,0.5*\linkpatterntangleheight);
\else
\iflinkpatternaxis
\draw (1,-0.5*\linkpatterntangleheight) -- (\lpsize,-0.5*\linkpatterntangleheight);
\draw (1,0.5*\linkpatterntangleheight) -- (\lpsize,0.5*\linkpatterntangleheight);
\else
\fi
\fi
\iflinkpatternrenumbered
\foreach\x/\xx in \linkpatternnumbering
{
\if\hasaprime\x %
\coordinate[vertex,lbla={\xx}] (v\x) at (\withoutprime{\x},0.5*\linkpatterntangleheight);
\else
\coordinate[vertex,lblb={\xx}] (v\x) at (\x,-0.5*\linkpatterntangleheight);
\fi
}
\else
\foreach\x in {1,...,\lpsize}
{
\coordinate[vertex,lbla={\bar\x}] (v\x') at (\x,0.5*\linkpatterntangleheight);
\coordinate[vertex,lblb={\x}] (v\x) at (\x,-0.5*\linkpatterntangleheight);
}
\fi
\foreach \a/\b in {#2}
{
  \gettikzxy{(v\a)}{\ax}{\ay}
  \gettikzxy{(v\b)}{\bx}{\by}
 \pgfmathsetmacro{\r}{(\ay+\by)/\linkpatternunit*\linkpatternlooseness*(-0.3+1/(1+abs(\ax-\bx)/\linkpatternunit))}
  \draw[edge] (v\a) .. controls (\ax,\r) and (\bx,\r) .. (v\b);
}
\iflinkpatterninverted
\end{scope}
\fi
\iflinkpatterntikzstarted\else%
\end{tikzpicture}%
\fi%
}}%
\newcommand\tanglelinkpattern[3][]{%
{%we want the macro modifications to be local only
\pgfkeys{/linkpattern/.cd,#1}
\iflinkpatterninverted%
\begin{tikzpicture}[/linkpattern/every linkpattern,yscale=-1]%
\else%
\begin{tikzpicture}[/linkpattern/every linkpattern]%
\fi%
\begin{scope}[yshift=-0.5*\linkpatternunit]
\tangle[#1,tikzstarted,halfnumbered,height=1]{#2}
\end{scope}
\linkpattern[#1,tikzstarted,numbered=false]{#3}
\end{tikzpicture}%
}}
\newcommand\circlelinkpattern[2][]{%in principle, could be merged
{%we want the macro modifications to be local only
\pgfkeys{/linkpattern/.cd,#1}%parse the list of optional arguments
\iflinkpatterntikzstarted\else%
\iflinkpatterninverted%
\begin{tikzpicture}[/linkpattern/every linkpattern,yscale=-1]%
\else%
\begin{tikzpicture}[/linkpattern/every linkpattern]%
\fi%
\fi%
\global\lpsize=\linkpatternsize
\edef\mylist{#2}%to avoid some annoying bugs
\foreach \x/\y in \mylist
{
\ifnum\x>\lpsize\global\lpsize=\x\fi
\ifnum\y>\lpsize\global\lpsize=\y\fi
}
\iflinkpatternaxis
\draw (0,0) circle (1);
\else\fi
\foreach\x in {1,...,\lpsize}
{
\pgfmathparse{(0.3*floor((\x-1)/\linkpatternfused)+0.7*((\x-0.5)/\linkpatternfused-0.5))*\linkpatternfused*360/\lpsize}
\coordinate[vertex] (v\x) at (\pgfmathresult:1);
}
\foreach \x/\y/\z in \mylist
{
\ifx\y\z%
\draw[edge] (v\x) .. controls ($0.5*(v\x)$) and  ($0.5*(v\y)$) .. (v\y);
\else
\draw[edge,decoration={markings,mark = at position 0.5 with { \arrow[semithick]{\z} }},postaction={decorate}] (v\x) .. controls ($0.5*(v\x)$) and  ($0.5*(v\y)$) .. (v\y);
\fi
}
\iflinkpatternnumbered%
\pgfmathparse{\lpsize/\linkpatternfused}
\lpsize=\pgfmathresult
\iflinkpatternrenumbered
\else
\def\linkpatternnumbering{1,...,\lpsize}
\fi
\newdimen\angle
\foreach\x/\xx in \linkpatternnumbering
{
  \pgfmathsetmacro{\angle}{360/\lpsize*(\x-1)}
  \node[outer sep=1pt,anchor=180+\angle] at (\angle:1) {$\scriptstyle\xx$}; % note the subtle anchoring: if we used for example "label", the anchoring would be restricted to the cardinal points i.e. north, north east, etc. instead we're using angle as anchor (documented in sec 39)
}
\else%
\fi%
\iflinkpatterntikzstarted\else%
\end{tikzpicture}%
\fi%
}}%
\makeatletter
\tikzset{circle split part fill/.style  args={#1,#2}{%
 alias=tmp@name,
  postaction={%
    insert path={
     \pgfextra{% 
     \pgfpointdiff{\pgfpointanchor{\pgf@node@name}{center}}%
                  {\pgfpointanchor{\pgf@node@name}{east}}%            
     \pgfmathsetmacro\insiderad{\pgf@x}
      \fill[#1] (\pgf@node@name.base) ([xshift=-\pgflinewidth]\pgf@node@name.east) arc
                          (0:180:\insiderad-\pgflinewidth)--cycle;
      \fill[#2] (\pgf@node@name.base) ([xshift=\pgflinewidth]\pgf@node@name.west)  arc
                           (180:360:\insiderad-\pgflinewidth)--cycle;                    }}}}}  
 \makeatother
\tikzset{bdot/.style={circle,circle split,draw,circle split part fill={black,white},thin,inner sep=1pt}}%
\tikzset{wdot/.style={circle,circle split,draw,circle split part fill={white,black},thin,inner sep=1pt}}%

\newdimen\unitsize
\newif\ifnumbered
\newcount\xx
\newcommand\rawpath[1]
{%
\xx=0
\coordinate[vertex] (current) at (0,0);
\foreach\x in {#1}
{
\global\advance\xx by 1
\if\x +%
\def\st{+1}
\else
\if\x 2%
\def\st{+1}
\else
\def\st{-1}
\fi
\fi
\ifnumbered
\draw[edge] (current) -- node[below,green!60!black] {$\scriptstyle \x$} ++(1,\st) coordinate[vertex] (current);
\else
\draw[edge] (current) -- ++(1,\st) coordinate[vertex] (current);
\fi
};
\ifnumbered
\draw[dotted] (0,0) -- (current |- 0,0);
\foreach\x in {1,...,\the\xx}
\path (\x,0)+(-0.5,0) node[below] {\x};
\fi
}
\newcommand\numberedpath[1]
{%
\begin{tikzpicture}[x=\unitsize,y=\unitsize,baseline=0]
\tikzset{vertex/.style={circle,fill=black,inner sep=1.2pt,outer sep=0pt}}%
\tikzset{edge/.style={very thick}}%
\numberedtrue\rawpath{#1}%
\useasboundingbox ([xshift=-2pt,yshift=-2pt] current bounding box.south west) 
([xshift=-2pt,yshift=2pt] current bounding box.north east);
\end{tikzpicture}%
}
\renewcommand\path[1]%this is a bit dangerous considering tikz also has a \path
{%
\begin{tikzpicture}[x=\unitsize,y=\unitsize,baseline=0]
\tikzset{vertex/.style={circle,fill=black,inner sep=1.2pt,outer sep=0pt}}%
\tikzset{edge/.style={very thick}}%
\numberedfalse\rawpath{#1}%
\useasboundingbox ([xshift=-2pt,yshift=-2pt] current bounding box.south west) 
([xshift=-2pt,yshift=2pt] current bounding box.north east);
\end{tikzpicture}%
}
\newcommand\treestyle{%
\tikzset{vertex/.style={circle,draw=black,thick,fill=green,inner sep=0pt,minimum size=4pt}}%
\tikzset{leaf/.style={regular polygon,regular polygon sides=3,draw=black,thick,fill=green,inner sep=0pt,minimum size=10pt}}%
\tikzset{root/.style={rectangle,draw=black,thick,fill=green,inner sep=0pt,minimum size=4pt,minimum width=12pt}}%
\tikzset{tree/.style={grow'=up,very thick,level distance=2.2\unitsize,sibling distance=2.2\unitsize}}%
}
\newdimen{\cellsize}
\newcommand\medboxes{\setlength{\cellsize}{14pt}\def\boxformat{}}

\medboxes
\tikzset{tableaubox/.style={draw=black,thin,solid,minimum size=\cellsize,inner sep=0pt}}
\def\cellbb{\useasboundingbox (-0.5*\cellsize,-0.5*\cellsize) rectangle ++(\cellsize,\cellsize);}
\tikzset{tableau/.style={matrix,name=tab,matrix anchor=tab-1-1.south west,inner sep=1pt,
execute at begin cell={\cellbb},execute at empty cell={\cellbb},matrix of math nodes,cells={anchor=center,draw=black,thin,solid,arrows=-},nodes={tableaubox,execute at begin node=\boxformat},nodes in empty cells}}
\def\activate#1{\begingroup
  \lccode`\~=`#1%
  \lowercase{\endgroup \let~#1}%
  \catcode`#1=13\relax}
\activate &
\def\tableau#1{\tikz[baseline=0]
\node[tableau]{#1};%
}

\begin{document}
\begin{abstract}
We provide simple rules for the computation of Kazhdan--Lusztig polynomials
in the maximal parabolic case. They are obtained by filling regions
delimited by paths with ``Dyck strips'' obeying certain rules. We compare
our results with those of Lascoux and Sch\"utzenberger.
\end{abstract}
\maketitle

\section{Introduction}

{\em Kazhdan--Lusztig polynomials}\/ 
were introduced in \cite{KL} as coefficients of the
change of basis from the standard basis of the Hecke algebra to a new
one, the Kazhdan--Lusztig basis. The latter was motivated by connections
to the representation theory of Weyl groups \cite{Spr76} %\cite{Jan79},%\cite{Spr78} 
and singularities of Schubert varieties \cite{KazLus80a} 
(see e.g. \cite{BilLak00} and references therein).
However, it reappeared since then in multiple contexts:
algebraic combinatorics \cite{KazLus80b} %,\cite{KasTan02}, 
Lie groups \cite{LusVog83}, %\cite{CasCol87},
the representation theory of Verma modules \cite{BryKas81,BelBer81},
and quantum groups \cite{Lus89}.   % \cite{KazLus91}
%modular representations \cite{LasLecThi96,Ari96}.

In \cite{Deo87} Deodhar introduced the concept of 
{\em parabolic Kazhdan--Lusztig polynomials}. % (see also \cite{Soe97}).
Roughly, they are associated to certain quotients of the regular representation
of the Hecke algebra ($q$-deformation of the induced representation
of one-dimensional representations of parabolic subgroups of the Coxeter group) 
in the same way as the usual Kazhdan--Lusztig polynomials
are associated to the regular representation, 
and the corresponding bases are projections of certain subsets
of the Kazhdan--Lusztig basis. Here we are concerned with type A
and a maximal parabolic subgroup, namely with Weyl group $\S_N$ and
the parabolic subgroup $\S_K\times \S_{N-K}$.
%Note that the basis in the type A parabolic case is also known to be
%the (dual) canonical basis of the quantum group in the sense of Lusztig 
%\cite{Lus90a,FreKhoKir98} 
%\cite{FreKho97}, \cite{FreKhoKir98}, 
%though this connection will not be described or used 
%in the present work.

The maximal parabolically induced representation of the Hecke algebra factors
through the {\em Temperley--Lieb algebra}~\cite{TL} and one 
expects simpler combinatorics than in the general case. 
Lascoux and Sch\"utzenberger~\cite{LS-KL} gave an algorithm to compute 
Kazhdan--Lusztig polynomials for Grassmannian permutations, which is 
equivalent to the maximal parabolic case (see \cite{Zel83} for a geometric 
interpretation). 
Also, there is a natural graphical description
of the basis and of the Temperley--Lieb action in terms of tangles
and link patterns, as used in models of two-dimensional statistical
mechanics \cite{BaxKelWu76,TL,Saleur-TL,Martin-TL}
and in knot theory~\cite{Kauf-TL}.
There is an abundant mathematical literature 
(see e.g. \cite{KL-KL,ES-KL,Bre94,Las95,Bre98,Bre02}) 
which provides explicit combinatorial formulae for some 
of these classes of Kazhdan--Lusztig polynomials. 

Due to the choice of the projection map %from $\R[\S_N]$ to $\R[\coset]$ 
(see Section \ref{secKL}), we have two types of parabolic Kazhdan--Lusztig 
polynomials studied in~\cite{LS-KL,Naruse01,Bre02}. 
The goal of the present paper is to provide a unified, self-contained 
treatment of maximal parabolic Kazhdan--Lusztig polynomials of both types
in the language of paths, 
similar to the one used by Brenti \cite{Bre02}. 
The main result is their computation according
to two graphical rules, denoted by I and II, rule II being equivalent
to Brenti's result.
The plan is as follows. In section 2, we introduce Kazhdan--Lusztig 
polynomials and their
maximal parabolic analogues and explain their duality.
Section 3 is the heart of the paper, in which 
we provide diagrammatic rules to compute the maximal parabolic 
polynomials.
In particular, the new rule (I) should be related to the one
given by Lascoux and Sch\"utzenberger in \cite{LS-KL}; 
and indeed, we provide a bijection between them in section 4.
We try to stay close to the conventions 
of the mathematical physics literature, 
which is where our motivation comes from.
More specifically, on the one hand
the Temperley--Lieb algebra and its ``link patterns'' have 
factorization properties
\cite{dG-review,MNdGB,KL-KL,dGP-factor}, which are relevant in calculations
that are performed in integrable loop models; on the other hand,
there are other explicit formulae \cite{artic41,artic42} which are
made in the ``standard basis'' of spin chains;
and we expect our formulae to be useful 
in connecting these different recent developments in integrable models. 

\section{Kazhdan--Lusztig polynomials} \label{secKL}
\subsection{Definition}
Given a positive integer $N$, 
we consider the symmetric group $\S_N$ with Coxeter generators
$s_i$, $i=1,\ldots,N-1$. Denote by $|v|$ the length of $v\in\S_N$,
that is the number of inversions $|v|:=\Card \{1\le i<j\le N, v(i)>v(j)\}$.
$\S_N$ is endowed with the (strong) {\em Bruhat order}
$\le$, that is $v\le w$ iff $v$ can be obtained by a series of
multiplications on the left (or right) by transpositions which each increase
length by one.

The {\em Hecke algebra}\/ $\H$ is the unital associative algebra over the ring $\R:=\mathbb{Z}[t,t^{-1}]$ 
with generators $T_i$, $i=1,\ldots,N-1$, and relations
\begin{align*}
 (T_i-t)(T_i+t^{-1})&=0 && 1\le i<N, \\
 T_iT_{i+1}T_i&=T_{i+1}T_iT_{i+1} && 1\le i<N-1,  \\
 T_iT_j&=T_jT_i && |i-j|>1.
\end{align*}

The standard basis $(T_v)_{v\in \S_N}$ of the Hecke algebra is obtained
by writing $T_v:=T_{i_1}\cdots T_{i_k}$ if $v=s_{i_1}\cdots s_{i_k}$ is
a reduced word in the elementary transpositions $s_i$ (see section 7 of \cite{Humph-refl}).

Define $a\to \overline a$ to be the involutive ring automorphism
of $\H$ such that $\overline{T_i}=T_i^{-1}$
and $\overline t=t^{-1}$ (that this map extends to a ring morphism
follows from invariance of the relations above by
$T_i\to T_i^{-1}$, $t\to t^{-1}$). Then 
\begin{theorem}[Kazhdan, Lusztig \cite{KL}] \label{kldef}
There exists a unique basis $(C_w)_{w\in\S_n}$ of $\H$  
such that $C_w=\overline{C_w}$ and
the matrix of change of basis $(P_{v,w})$
from the $T_v$ to the $C_w$ is ``upper triangular'' w.r.t.\ Bruhat order, i.e.\ 
\[
C_w=\sum_{\substack{v\in \S_N\\ v\le w}} P_{v,w}(t^{-1}) T_{v}
\]
where the polynomials $P_{v,w}(t^{-1})\in t^{-1}\Z[t^{-1}]$ if $v<w$ and $P_{v,v}=1$.
\end{theorem}
In fact, $\deg P_{v,w}\le |w|-|v|$, and
the Kazhdan--Lusztig (KL) polynomials are by definition the polynomials
$t^{|w|-|v|} P_{v,w}(t^{-1})\in \Z[t^2]$. 

\rem{we got different convention than KL. their $q$ is like our
$t^{-2}$ (hence their obnoxious use of the bar in the definition). 
the other difference is due to their normalization of
Hecke generators.}

\subsection{Maximal parabolic case}
Given $0\le K\le N$, 
we now consider the subgroup
$\S_K\times \S_{N-K}\subset \S_N$ 
with generators $s_i$, $i=1,\ldots,K-1,K+1,\ldots,N-1$.
The set of left cosets $\coset$ has a natural induced order:
$x\le y$ iff there exist $v\in x,w\in y$ such that $v\le w$, and a length:
$|x|=\min_{v\in x} |v|$.

Let us define $\M$ to be a free $\R$-module with basis
indexed by $\coset$: $\M:=\langle m_{x}, x\in\coset\rangle$.
The projection $\varphi$ from $\S_N$ to $\coset$ induces 
{\em two}\/ natural projection maps $\varphi^{\pm}$
from 
%$\H \cong \R[\S_N]$ 
$\H$
to 
%$\M^\pm\cong\R[\coset]$
$\M$,
given by $\varphi^{\pm}(T_v):=(\pm t^{\pm 1})^{|v|-|\varphi(v)|} m_{\varphi(v)}$. 
Fix $\epsilon\in\{+,-\}$.
In order to define a representation of $\H$ on the $R$-module $\M$,
we require that $\varphi^\epsilon$
commutes with the action of the Hecke algebra
(cf lemma 2.2 of \cite{Deo87}), where the latter acts on itself by left
multiplication; this leads to the following action
of the generators $T_i$ on $\M$:
\begin{equation}\label{heckeaction}
	T_i m_{x}  
	=
	\begin{cases}
	 \epsilon t^{\epsilon} m_x  & s_ix=x, \\
	 m_{s_i x} & s_i x\ne x, |s_i x|>|x|\\
	 (t-t^{-1})m_{x} + m_{s_i x} & s_i x\ne x, |s_i x|<|x|.\\
         \end{cases}
\end{equation}
This endows $\M$ with the structure of an $\H$-module, which
is denoted by $\M^\epsilon$.
\rem{+ is LS, - is Brenti, which is why $t=-q^{-1}$}

Similarly, requiring that $\varphi^\epsilon$ commute with the bar involution 
defines uniquely its action on $\M^\epsilon$. 

We can now define parabolic analogues of KL basis and polynomials:
\begin{theorem}[Deodhar \cite{Deo87}]\label{klparadef}
There exists a unique basis $(C^\pm_x)_{x\in\coset}$ of $\M^\pm$  
such that $C^\pm_x=\overline{C^\pm_x}$ and
the matrix of change of basis $(P^\pm_{x,y})$
from the $m_x$ to the $C^\pm_y$ is ``upper triangular'', i.e.
\[
C^{\pm}_{y} = \sum_{\substack{x\in\coset\\x\le y}} P^{\pm}_{x,y}(t^{-1}) m_{x}
\]
where the polynomials $P^\pm_{x,y}(t^{-1})\in t^{-1}\Z[t^{-1}]$ 
if $x<y$ and $P_{x,x}=1$.
\end{theorem}
In fact, $\deg P^\pm_{x,y}\le |y|-|x|$, and
the parabolic Kazhdan--Lusztig polynomials are 
by definition the %polynomials
$t^{|y|-|x|} P^{\pm}_{x,y}(t^{-1})\in \Z[t^2]$. 
Here we prefer to use directly the polynomials $P^{\pm}_{x,y}(t^{-1})$.

\subsection{Combinatorial description}\label{secbij}
There are various ways to describe explicitly the cosets
in $\coset$. We are of course mostly interested in their {\em path}\/
representation, but we discuss in this section other useful descriptions.

Let $\epsilon\in \{-,+\}$. We consider the following sets:
\begin{enumerate}[(1)]
\setcounter{enumi}{-1}
\item $\coset$ \\ and the sets of:
\item {\bf Binary strings}, i.e.\ elements of $\{1,2\}^N$,
such that there are $K$ $1$'s and $N-K$ $2$'s.
\item {\bf Paths} from $(0,0)$ to $(N,\epsilon(2K-N))$ with steps
$(1,\pm 1)$.
\item {\bf Ferrers diagrams} inside the rectangle $K\times(N-K)$.
\item {\bf Link patterns} with at most $\min(K,N-K)$ pairings, 
where link patterns are planar pairings of a subset
of $\{1,\ldots,N\}$ in such a way that unpaired vertices belong
to the infinite connected component. %(see the example below).
\item {\bf (anti)Grassmannian permutations}, that is permutations
$\sigma$ such that
$1\le i<j\le K$ or $K+1\le i<j\le N$ implies $\epsilon\sigma(i)>
\epsilon\sigma(j)$.
\item {\bf Standard Young tableaux with at most two rows} (resp.\ two columns
for $\epsilon=+$)
and with $N$ boxes, 
whose second row (resp.\ column) is of length less or equal to $\min(K,N-K)$.
\end{enumerate}

as well as the following maps between these sets:
\begin{itemize}
\item[(0)$\to$(1):] such binary strings are the orbits under the natural action 
of $\S_N$ on $\{1,2\}^N$, with representative
$(\underbrace{1,\cdots,1}_K,\underbrace{2,\cdots,2}_{N-K})$.
The latter has stabilizer $\S_K\times \S_{N-K}$.

\item[(1)$\to$(2):] a sequence $v\in\{1,2\}^N$ is identified
with the path with $i^{\rm th}$ step $(1,\epsilon(-1)^{1+v_i})$.

\item[(2)$\to$(3):] to a path is associated the (45 degrees rotated) Ferrers diagram
located between it and the smallest path for $\le$ (corresponding to the binary
string $(1,\ldots,1,2,\ldots,2)$ and to the coset of the identity;
it is the lowest path for $\epsilon=-$, the highest path for $\epsilon=+$).

\item[(2)$\to$(4):]
pair midpoints of steps of equal height
such that the horizontal segment between them stays strictly below
the path. (see the example below)

\item[(0)$\to$(5):] in each coset $x$, there is exactly one 
Grassmannian permutation,
denoted by $\short{x}$: 
it is the ``shortest representative'' (of shortest length).
Note that by definition, $|x|=|\short{x}|$, 
and $x\le y$ iff $\short{x}\le\short{y}$. 
%By abuse of notation
%we shall also abbreviate $\short{u}:=\short{\varphi(u)}$ for $u\in\S_N$.
Inversely there is exactly one anti-Grassmannian permutation in each coset:
it is the ``longest representative'',
and can be written $\short{x}\tilde w_0$, where $\tilde w_0$
is the longest element of $\S_K\times \S_{N-K}$, namely
$\left(\begin{smallmatrix}
K&\cdots&1&N&\cdots&K+1\\
1&\cdots&K&K+1&\cdots&N
\end{smallmatrix}\right)$.

\item[(5)$\to$(6):] 
applying the Robinson--Schensted algorithm to $\short{x}$
results in a pair of Young tableaux of same shape with at most
2 rows; keep only the first tableau, the second one being
entirely fixed by its shape, say $(N-i,i)$, to be:
\setlength{\cellsize}{19pt}\def\boxformat{\scriptstyle}
$\tableau{1&2&\cdots&K&
\node[align=center]{$\scriptstyle K+i$\\[-7pt]$\scriptstyle +1$};
&\cdots&N\\ K+1&\cdots &K+i\\}$.
Similarly,
applying the Robinson--Schensted algorithm to $\short{x}\tilde w_0$
results in a pair of Young tableaux of same shape with at most
2 columns; keep only the first tableau, the second one being
entirely fixed by its shape. Note that these two tableaux
are {\em not}\/ transpose of each other.
\end{itemize}

\begin{lemma}\label{combidescr}
The maps described above are bijections.
\end{lemma}
The proofs are standard (see e.g.\ \cite{Stan-EC2} and online
supplements at \cite{Stan-ECsup}).

\begin{example} 
We choose the $\epsilon=-$ convention to draw (2), (4).
%\numberedpath{+,-,-,+,+,-,+,-,-,-}

\noindent (2) $\leftrightarrow$ (1):
\numberedpath{2,1,1,2,2,1,2,1,1,1}
\quad$\leftrightarrow$\quad
(2,1,1,2,2,1,2,1,1,1)

%\item[(2) $\leftrightarrow$ (3)]
\noindent (2) $\leftrightarrow$ (3)
\hbox{\unitsize=0.5cm%
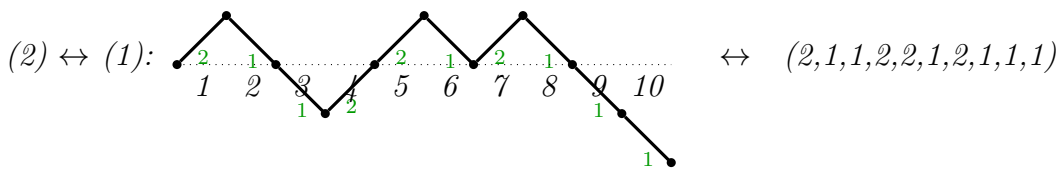
\begin{tikzpicture}[x=\unitsize,y=\unitsize,baseline=0]
\tikzset{vertex/.style={circle,fill=black,inner sep=1.2pt,outer sep=0pt}}%
\tikzset{edge/.style={very thick}}%
\numberedfalse%
\rawpath{+,-,-,+,+,-,+,-,-,-}%
\rawpath{-,-,-,-,-,-,+,+,+,+}%
\path[use as bounding box] (-2pt,-2\unitsize) rectangle (10\unitsize+2pt,\unitsize);
\draw[dashed] (1,-1) -- (2,0) (2,-2) -- (3,-1) (3,-3) -- (6,0) (4,-4) -- (8,0) (5,-5) -- (9,-1);
\draw[dashed] (3,-1) -- (7,-5) (4,0) -- (8,-4) (6,0) -- (9,-3);
\end{tikzpicture}%
\qquad$\leftrightarrow$\qquad\tableau{&&&&&&\\&&&\\&&&\\&&\\}}

\noindent (2) $\leftrightarrow$ (4):
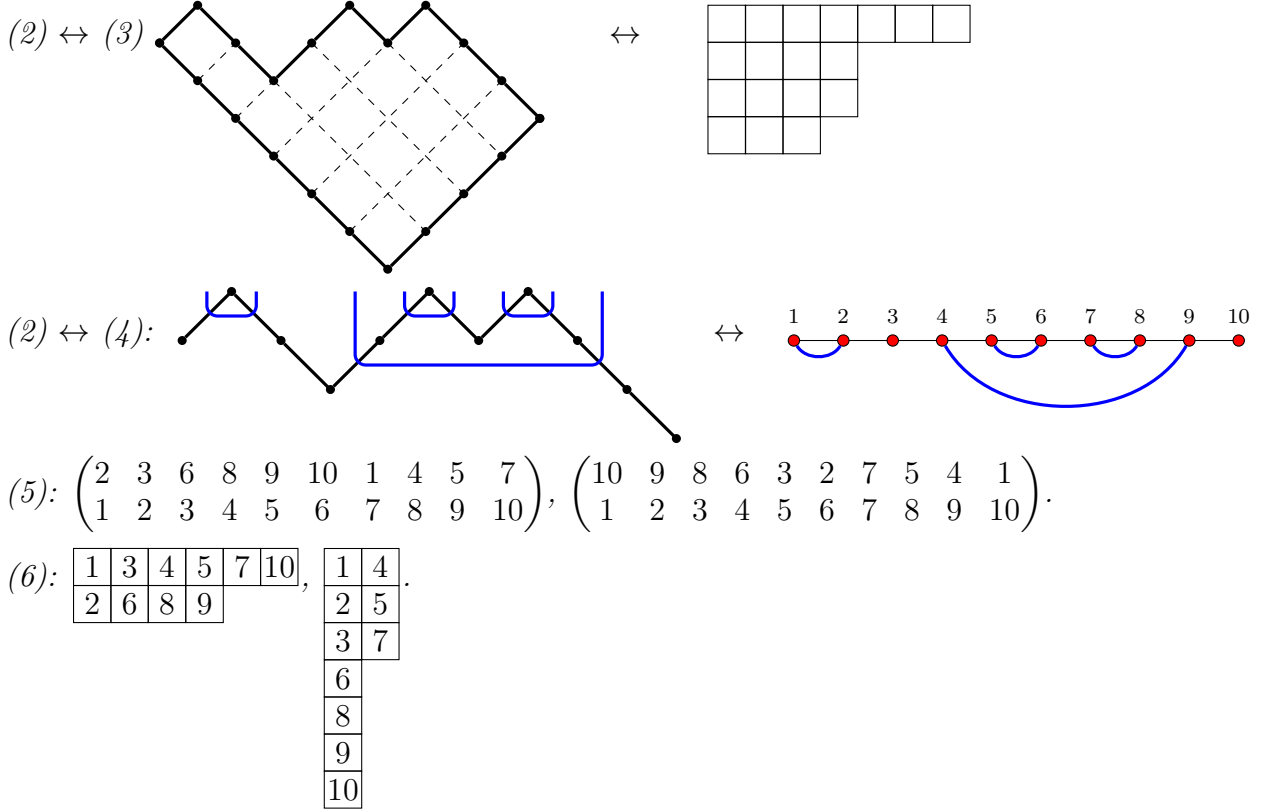
\hbox{
\begin{tikzpicture}[x=\unitsize,y=\unitsize,baseline=0]
\tikzset{vertex/.style={circle,fill=black,inner sep=1.2pt,outer sep=0pt}}%
\tikzset{edge/.style={very thick}}%
\numberedfalse%
\rawpath{+,-,-,+,+,-,+,-,-,-}%
\draw[rounded corners,blue,very thick] (0.5,1) -- (0.5,0.5) -- (1.5,0.5) -- (1.5,1);
\draw[rounded corners,blue,very thick] (4.5,1) -- (4.5,0.5) -- (5.5,0.5) -- (5.5,1);
\draw[rounded corners,blue,very thick] (6.5,1) -- (6.5,0.5) -- (7.5,0.5) -- (7.5,1);
\draw[rounded corners,blue,very thick] (3.5,1) -- (3.5,-0.5) -- (8.5,-0.5) -- (8.5,1);
\end{tikzpicture}%
\quad$\leftrightarrow$\quad\linkpattern[unit=\unitsize,numbered,inverted,unpaired length=0]{1/2,4/9,5/6,7/8,10/10}}

\noindent (5): $\left(
\begin{matrix}
2&3&6&8&9&10&1&4&5&7\\
1&2&3&4&5&6&7&8&9&10
\end{matrix}
\right)$,
$\left(\begin{matrix}
10&9&8&6&3&2&7&5&4&1\\
1&2&3&4&5&6&7&8&9&10
\end{matrix}\right)$.

\noindent (6): $\medboxes\tableau{1&3&4&5&7&10\\2&6&8&9\\}$, 
$\medboxes\tableau{1&4\\2&5\\3&7\\6\\8\\9\\10\\}$.
\end{example}

In what follows, we shall mostly use
the path representation, or interchangeably 
the closely related Ferrers diagram representation. More precisely
the bijection to paths with the sign convention $\epsilon\in\{-,+\}$
will be used to index bases of $\M^{\epsilon}$. 
The set of paths from $(0,0)$ to $(N,2K-N)$ will be denoted by $\mathcal{P}_{N,K}$.
%In order to avoid
%any confusion, greek letters will be used to denote paths.

It is perhaps useful to rewrite the action \eqref{heckeaction} of the Hecke algebra
on $\M^\pm$ in terms of local changes of paths: (only the steps $i$ and $i+1$ are depicted)
{\setlength{\unitsize}{0.2cm}
\def\moddots#1{\raise#1\unitsize\hbox{$\scriptstyle\ldots$}}
\def\pp{{\ldots\path{+,+}\moddots2}}
\def\mm{{\ldots\path{-,-}\moddots{-2}}}
\def\mp{{\ldots\path{-,+}\ldots}}
\def\pm{{\ldots\path{+,-}\ldots}}
\begin{align*}
\epsilon=-: \qquad
T_i\ m_\pp &= -t^{-1}\ m_\pp \\
T_i\ m_\mm &= -t^{-1}\ m_\mm \\
T_i\ m_\mp&=  m_\pm \\
T_i\ m_\pm&=  (t-t^{-1})\ m_\pm
+m_\mp \\[4pt]
\epsilon=+: \qquad
T_i\ m_\pp &= t\ m_\pp \\
T_i\ m_\mm &= t\ m_\mm \\
T_i\ m_\pm&=  m_\mp \\
T_i\ m_\mp&=  (t-t^{-1})\ m_\mp
+m_\pm \\
\end{align*}
}
In terms of the associated Ferrers diagrams, the third and fourth lines involves
adding and removing a box, respectively.

We have the following additional easy lemma:
\begin{lemma}\label{easylemma}
Let $x,y\in\coset$ and call $\alpha,\beta$ the associated paths
with convention $\epsilon$.
Then $x\le y$ iff $\alpha$ is {\em below} $\beta$
for $\epsilon=-$, above $\beta$ for $\epsilon=+$; and $|x|$ is the
number of boxes of the corresponding Ferrers diagram, also denoted by
$|\alpha|$ in what follows.
\end{lemma}
\rem{\begin{figure}[h!]
\[
\begin{tikzpicture}[x=\unitsize,y=\unitsize,baseline=0]
\tikzset{vertex/.style={circle,fill=black,inner sep=1.2pt,outer sep=0pt}}%
\tikzset{edge/.style={very thick}}%
\numberedfalse%
\rawpath{+,-,-,+,+,-,+,-,-,-}%
\rawpath{+,-,-,+,-,-,+,-,+,-}%
\numberedfalse%
\tikzset{vertex/.style={}}%
\tikzset{edge/.style={dashed}}%
\rawpath{+,-,-,+,-,+,-,+,-,-}
\end{tikzpicture}
\]
\caption{An example of one path above another path, with length difference 4.
\rem{is this figure really useful?}}
\label{figtwopaths}
\end{figure}}
\begin{proof}
Let us prove the case $\epsilon=-$.
If the path $\beta$ is above $\alpha$ then $y$ can be
obtained from $x$ by a series of multiplications on the left by elementary
transpositions (as mentioned above, this corresponds to adding one box at a time
on top of the path). %, see also Fig.~\ref{figtwopaths}
Therefore $x\le y$.
Conversely, assume $x\le y$, i.e., $u:=\short{x}\le v:=\short{y}$. 
Define the height function associated to a permutation by
$h(x)_{i,j}=\#\{k\le i: x(k)\le j\}$,
$i,j=0,\ldots,N$.
Then it is well-known that $u\le v$ iff $h(u)_{i,j}\ge h(v)_{i,j}$ for all
$i,j=0,\ldots,N$. And the path associated to a Grassmannian permutation
$u$ is nothing but the path with set of vertices
$\{(j,j-2h(u)_{K,j}), j=0,\ldots,N\}$. This proves the assertion.

The first part of the reasoning also shows that increasing the length by one
is the same as adding one box under the path, which leads to the second
part of the lemma. 
\end{proof}

\subsection{Connection between KL and parabolic KL polynomials}
Since the projections $\varphi^{\pm}$ commute with the Hecke action
and with the bar involution,
images of elements $C_w$ of the Kazhdan--Lusztig basis of $\H$
under $\varphi^{\pm}$ are natural candidates for their parabolic
counterparts $C_w^{\pm}$. And indeed, one can show that
$\varphi^\epsilon(C_w)=C^\epsilon_{\varphi(w)}$ 
if $w$ is the shortest (Grassmannian) representative of its coset
for $\epsilon=-$, and the longest representative for 
$\epsilon=+$.
Note however that the definitions in Theorems \ref{kldef} and
\ref{klparadef} of KL bases 
break the symmetry in the definition of the Hecke algebra
between $t$ and $-t^{-1}$ 
(by requiring the coefficients to be polynomials in $t^{-1}$)
which is therefore not apparent in the resulting formulae
for parabolic KL polynomials:
\begin{prop}[Deodhar \cite{Deo87}]\label{kltoparakl}
\begin{align*}
P^{+}_{x,y}&=P_{v,w} && \text{$v=\short{x}\tilde w_0$, $w=\short{y}\tilde w_0$ longest representatives}\\
P^{-}_{x,y}&=\sum_{v\in x} (-t)^{|x|-|v|} P_{v,w}
&& \text{$w=\short{y}$ shortest representative}\end{align*}
\end{prop}

\subsection{Duality}
There is a general duality satisfied by KL polynomials
(Theorem 3.1 of \cite{KL}). Let $w^0$ be the longest element of $\S_N$,
namely $w^0=\left(\begin{smallmatrix}N&\cdots&1\\ 1&\cdots&N\end{smallmatrix}\right)$. Reformulated in our language, this result becomes
\begin{theorem}[Kazhdan, Lusztig \cite{KL}]\label{klduality}
The following inversion formulae hold:
\begin{align*}
\sum_{w\in S_N} (-1)^{|v|+|w|} P_{u,w} P_{w^0 v,w^0 w}&=
\delta_{u,v}\\
\sum_{u\in S_N} (-1)^{|u|+|v|} P_{u,v} P_{w^0 u,w^0 w}&=
\delta_{v,w}
\end{align*}
\end{theorem}

For our purposes it is more convenient to have $w^0$ act on the right,
which amounts to using the opposite product, or to applying the small
\begin{lemma}\label{sharpinv}
Let $u^\sharp = w^0 u w^0$. Then
\[
P_{u^\sharp,v^\sharp}=P_{u,v}
\]
\end{lemma}
\begin{proof}
Firstly, $\sharp$ preserves the Bruhat order.
Secondly, extend $\sharp$ into an involution of $\H$ with
$T_v^\sharp=T_{v^\sharp}$. Noting that $\sharp$ and bar involutions
commute, we conclude that $C_{w^\sharp}=(C_w)^\sharp$,
hence the result.
\end{proof}

Recall that we also have the longest element
in $\S_K\times \S_{N-K}$:
$\tilde w_0=\left(\begin{smallmatrix}
K&\cdots&1&N&\cdots&K+1\\
1&\cdots&K&K+1&\cdots&N
\end{smallmatrix}\right)$. Write $w^0=\eta \tilde w^0$,
where $\eta=\left(\begin{smallmatrix}
N-K+1&\cdots&N&1&\cdots&N-K\\
1&\cdots&K&K+1&\cdots&N
\end{smallmatrix}\right)$.

We now switch as promised to the path indexation. All the paths
in this section are in $\mathcal{P}_{N,K}$, i.e., from $(0,0)$ to $(N,2K-N)$. 
Let $\gamma$ be such a path.
According to lemma \ref{combidescr},
they can be interpreted as either
($\epsilon=-$) a shortest representative in $\cosetb$, say $w$, 
or ($\epsilon=+$) a longest representative in $\coset$, 
say $w'$. The key remark is that we have $w'=ww_0$: indeed
multiplying by $\eta$ on the right
flips the path upside down (following the different convention for
paths depending on $\epsilon$), and multiplying
by $\tilde w_0$ turns shortest into longest representative.
Therefore, given two paths $\beta,\gamma\in\mathcal{P}_{N,K}$, 
one can associate to them
$v$ and $w$, the shortest representatives as above, 
and write, applying Proposition \ref{kltoparakl}:
\begin{align*}
\sum_{\alpha\in\mathcal{P}_{N,K}} (-1)^{|\alpha|+|\beta|}
P^-_{\alpha,\beta} P^+_{\alpha,\gamma} 
&=\sum_{z\in \coset} \sum_{u\in z} (-1)^{|u|+|v|} 
P_{u,v}\, t^{|z|-|u|} P_{\short{z} w^0,w w^0}
\\
&=\sum_{z\in \coset} \sum_{u\in z} (-1)^{|u|+|v|} 
P_{u,v} P_{u w^0,w w^0}
&&\text{by (2.3.g) of \cite{KL}}
\\
&=\sum_{u\in\S_N} (-1)^{|u|+|v|}
P_{u,v} P_{u w^0,w w^0}
\end{align*}
where in the application of (2.3.g) of \cite{KL} we set
$x=u w^0$, $y=w w^0$ and use the opposite product.

Writing that $P_{u w^0,w w^0}=P_{w^0 u,w^0 w}$ (lemma \ref{sharpinv})
leads to the second identity of theorem \ref{klduality}, so that
\[
\sum_{\alpha\in\mathcal{P}_{N,K}} (-1)^{|\alpha|+|\beta|}
P^-_{\alpha,\beta} P^+_{\alpha,\gamma} 
=\delta_{\beta,\gamma}
\]

Note that $(-1)^{|\alpha|+|\beta|} P^-_{\alpha,\beta}(t^{-1})=P^-_{\alpha,\beta}(-t^{-1})$. We reach the result
\begin{theorem}\label{paraduality}
\[
\sum_{\alpha\in\mathcal{P}_{N,K}}
P^-_{\alpha,\beta}(-t^{-1}) P^+_{\alpha,\gamma}(t^{-1}) 
=\delta_{\beta,\gamma}
\]
\end{theorem}

In other words, the change of basis for $\M^+$ is up to $t\to -t$ the inverse
transpose of the one for $\Md^-$, which is just
a manifestation of usual (linear algebra) duality.

\section{Path representation}

\subsection{Dyck strips}

\begin{figure}[h!]
\setlength{\unitsize}{0.35cm}%
\tikzset{vertex/.style={}}%
\begin{tikzpicture}[x=\unitsize,y=\unitsize,baseline=0]
\tikzset{edge/.style={}}%
\rawpath{+,-}\rawpath{-,+}
\end{tikzpicture}
\qquad
\begin{tikzpicture}[x=\unitsize,y=\unitsize,baseline=0]
\tikzset{edge/.style={}}%
\rawpath{+,+,+,-,-,-}\rawpath{-,+,+,-,-,+}
\tikzset{edge/.style={dotted}}%
\rawpath{+,-,+,+,-,-}\rawpath{+,+,-,-,+,-}
\end{tikzpicture}
\qquad
\begin{tikzpicture}[x=\unitsize,y=\unitsize,baseline=0]
\tikzset{edge/.style={}}%
\rawpath{+,+,-,+,-,-}\rawpath{-,+,-,+,-,+}
\tikzset{edge/.style={dotted}}%
\rawpath{+,-,+,-,+,-}
\end{tikzpicture}
\qquad
\begin{tikzpicture}[x=\unitsize,y=\unitsize,baseline=0]
\tikzset{edge/.style={}}%
\rawpath{+,+,+,-,+,-,-,-}\rawpath{-,+,+,-,+,-,-,+}
\tikzset{edge/.style={dotted}}%
\rawpath{+,+,-,+,-,+,-,-}\rawpath{+,-,+,-,+,-,+,-}
\end{tikzpicture}
\caption{Some Dyck strips.}\label{fig-d-strip}
\end{figure}
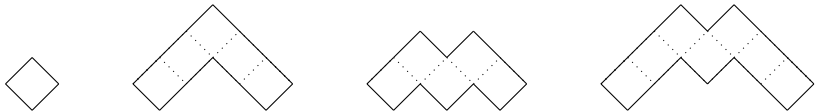

{\it A Dyck path}\/ of length $2l, l\ge0$,
is a path from some 
$(x,y)\in\mathbb{Z}^2$ to $(x+2l,y)$ and not crossing below the 
horizontal line at height $y$. 
{\it A Dyck strip}\/ of length $2l+1$ is obtained 
by putting unit boxes (45 degrees 
rotated) whose centers are at the vertices of a Dyck path of length $2l$ (see 
some examples on Fig~\ref{fig-d-strip}). 

Hereafter, a box $(x,y)$ means a unit box whose center is $(x,y)$. 
Let $b$ be a box $(x,y)$. 
Four boxes $(x\pm1,y\pm1)$ are neighbors of $b$.  
The box $(x+1,y+1)$ is said to be NE (north-east) of $b$ 
and similarly the other three boxes are NW, SW and SE of $b$.
The two boxes $(x,y\pm2)$ are said to be {\it just above}\/ or 
{\it just below}\/ $b$.

We now define two relations on the set of Dyck strips
as follows. Given an ordered pair of such Dyck strips $(D,D')$, 
we say that they satisfy rule I/II iff:
\begin{description}
 \item[Rule I]  
 If there exists a box of ${D}$ just below a box of ${D'}$, then all boxes just below 
a box of ${D'}$ belong to ${D}$.
\item[Rule II] 
If there exists a box of ${D'}$ just above, NW or NE of a box of ${D}$, then all boxes 
just above, NW and NE of a box of ${D}$ belong to ${D}$ or ${D'}$.
\end{description}
We are interested in decomposing skews into unions of strips
according to one of these rules.
Roughly, Rule I (resp.\ Rule II) means that we are allowed to pile 
Dyck strips of smaller or equal (resp.\ longer) length on top of a Dyck strip. 

\begin{figure}
\setlength{\unitsize}{0.25cm}%
\tikzset{vertex/.style={}}%
\tikzset{edge/.style={}}%
\begin{tikzpicture}[x=\unitsize,y=\unitsize,baseline=0]
\rawpath{+,+,+,+,+,+,+,+,+,-,-,-,+,-,-,-,-,-,-}
\rawpath{+,+,+,+,+,+,+,-,+,-,+,-,+,-,-,-,-,-,-}
\rawpath{+,+,+,+,+,+,-,+,+,-,-,+,+,-,-,-,-,-,-}
\rawpath{+,+,+,+,+,-,+,-,+,-,+,-,+,-,+,-,-,-,-}
\rawpath{+,+,+,+,-,+,-,+,+,-,+,-,+,-,-,+,-,-,-}
\rawpath{+,+,+,+,-,-,+,+,+,-,+,-,+,-,-,-,+,-,-}
\rawpath{+,+,+,+,-,-,-,+,+,-,+,-,+,-,-,+,-,+,-}
\rawpath{+,-,+,+,-,-,+,+,+,-,+,-,+,-,-,-,+,-,+}
\rawpath{-,+,+,+,-,-,+,+,+,-,+,-,+,-,-,-,-,+,+}
\rawpath{-,-,+,+,-,-,+,+,+,-,+,-,+,-,-,-,+,+,+}
\end{tikzpicture}
\qquad
\begin{tikzpicture}[x=\unitsize,y=\unitsize,baseline=0]
\rawpath{+,+,+,-,+,-,-,+,+,-,-,-}
\rawpath{-,+,+,-,+,-,-,+,+,-,-,+}
\rawpath{-,-,+,-,+,-,+,-,+,-,+,+}
\rawpath{-,-,-,+,+,-,+,-,+,-,+,+}
\end{tikzpicture}
\caption{Examples of stacks of Dyck strips satisfying rule I (left) and rule II (right).}\label{fig-x}
\end{figure}
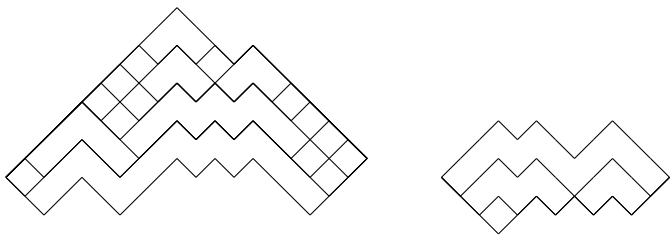

\newcommand\CI{\mathrm{Conf}^{I}}
\newcommand\CII{\mathrm{Conf}^{II}}

Let $\alpha,\beta\in\mathcal{P}_{N,K}$ be two paths as defined in
section~\ref{secbij}. 
We consider filling the closed domain between these two paths with Dyck strips
(such that no Dyck strips overlap, and every unit box is filled). 
Let us denote by $\mathrm{Conf}(\alpha,\beta)$ the set of all 
such possible configurations of Dyck strips, 
and $\mathrm{Conf}^{I/II}(\alpha,\beta)$ the subset of configurations satisfying rule I/II.
We denote the number of Dyck strips in a configuration $\mathcal{D}$
by $|\mathcal{D}|$.

\begin{defn}
The generating function of Dyck strips for the paths $\alpha<\beta$ 
in $\mathcal{P}_{N,K}$ is defined as 
\[
 Q_{\alpha,\beta}^{X,\epsilon}(t^{-1})=\sum_{\mathcal{D}\in\mathrm{Conf}^{X}(\alpha,\beta)}
t^{-|\mathcal{D}|}.
\]
where $X=I,II$ and $\epsilon=\pm$ is the order convention 
as in Lemma \ref{easylemma}.
$Q_{\alpha,\alpha}^{X,\epsilon}(t^{-1})=1$ and $Q_{\alpha,\beta}(t^{-1})=0$ if
$\alpha\not\le\beta$.
\end{defn}
%

%Note that in the case of 
Due to the implied condition of Rule II, we have at most one configuration 
for given paths $\alpha$ and $\beta$ due to the condition $l(D)<l(D')$. 
In other words, the $Q^{II}$ are monomials.
%

%Once we fix a path in the $K\times(N-K)$ rectangle, the choice of 
%$\epsilon$ changes the corresponding representative (shortest or longest) in 
%$\mathcal{S}_N/(\mathcal{S}_N\times\mathcal{S}_{N-K})$ and equivalently changes% the pictorial 
%order of paths (above or below) corresponding to the Bruhat order. 
Recall that according to Lemma \ref{easylemma},
$\alpha\le\beta$ means $\alpha$ is pictorially above 
(resp.\ below) $\beta$ for $\epsilon=+$ (resp. $\epsilon=-$). 
Therefore, it is obvious that 
\begin{lemma}\label{lemma-flip}
\[
 Q^{X,+}_{\alpha,\beta}(t^{-1})= Q^{X,-}_{\beta,\alpha}(t^{-1}). 
\]
\end{lemma}

\begin{example}
When $(\alpha,\beta)=(111212222,211212221)$,
\[
\CI(\alpha,\beta)=
\left\{
\setlength{\unitsize}{0.2 cm}%
\tikzset{vertex/.style={}}%
\tikzset{edge/.style={}}%
\begin{tikzpicture}[x=\unitsize,y=\unitsize,baseline=0]
\rawpath{+,+,+,-,+,-,-,-,-}
\rawpath{+,+,-,+,-,+,-,-,-}
\rawpath{+,-,+,-,+,-,+,-,-}
\rawpath{+,-,+,-,+,-,-,+,-}
\rawpath{-,+,+,-,+,-,-,-,+}
\end{tikzpicture}
,
\begin{tikzpicture}[x=\unitsize,y=\unitsize,baseline=0]
\rawpath{+,+,+,-,+,-,-,-,-}
\rawpath{+,+,+,-,-,+,-,-,-}
\rawpath{+,-,+,-,+,-,+,-,-}
\rawpath{+,-,+,-,+,-,-,+,-}
\rawpath{-,+,+,-,+,-,-,-,+}
\end{tikzpicture}
,
\begin{tikzpicture}[x=\unitsize,y=\unitsize,baseline=0]
\rawpath{+,+,+,-,+,-,-,-,-}
\rawpath{+,+,-,+,+,-,-,-,-}
\rawpath{+,-,+,-,+,-,+,-,-}
\rawpath{+,-,+,-,+,-,-,+,-}
\rawpath{-,+,+,-,+,-,-,-,+}
\end{tikzpicture}
,
\begin{tikzpicture}[x=\unitsize,y=\unitsize,baseline=0]
\rawpath{+,+,+,-,+,-,-,-,-}
\rawpath{+,-,+,-,+,-,+,-,-}
\rawpath{-,+,+,-,+,-,-,+,-}
\rawpath{-,+,+,-,+,-,-,-,+}
\end{tikzpicture}
,
\begin{tikzpicture}[x=\unitsize,y=\unitsize,baseline=0]
\rawpath{+,+,+,-,+,-,-,-,-}
\rawpath{-,+,+,-,+,-,-,+,-}
\rawpath{-,+,+,-,+,-,-,-,+}
\end{tikzpicture}
\right\}
\]
The corresponding
generating function is $Q^{I,+}_{\alpha,\beta}(t^{-1})=t^{-8}(1+2t^2+t^4+t^6)$.
%,whereas $Q^{II,+}_{\alpha,\beta}(t^{-1})=0$.
\end{example}

%\subsection{KL basis}

The relations among the Kazhdan--Lusztig polynomials $P^{\pm}_{\alpha,\beta}$ and 
the generating functions $Q^{X,\epsilon}_{\alpha,\beta}$ that we shall establish 
in subsequent sections are summarized as:
\[
 \begin{CD} \begin{array}{c}
 \mbox{Theorem~\ref{cor-minus-KL}} \\
 P^-_{\alpha,\beta}=Q^{II,-}_{\alpha,\beta} 
 \end{array} 
@.  \xlongleftrightarrow{\quad\mbox{transpose}\quad} @.
\begin{array}{c}
 \mbox{Corollary~\ref{inverse-LS}} \\
Q^{II,+}_{\alpha,\beta} 
\end{array} \\ 
 \mbox{$\Bigg\updownarrow$inverse} @. @. \mbox{$\Bigg\updownarrow$inverse}  \\
\begin{array}{c}
 \mbox{Theorem~\ref{mainthm}} \\
Q^{I,-}_{\alpha,\beta} \end{array} @.
 \xlongleftrightarrow{\quad\mbox{transpose}\quad} @.
\begin{array}{c}
 \mbox{Corollary~\ref{theorem-LS} } \\
P^+_{\alpha,\beta}=Q^{I,+}_{\alpha,\beta} 
\end{array}
 \end{CD}
\]

\subsection{On the module $\mathcal{M}^-_{N,K}$}
For the purposes of this section, 
we identify a path and binary string of $1$ and $2$ with convention
$-$ as in section \ref{secbij}. 
We denote by $\alpha=\alpha_1\ldots\alpha_N$ a binary string of $N$ letters. 
\begin{defn}
For given paths $\alpha,\beta\in \mathcal{P}_{N,K}$, 
we define 
\[
 d(\alpha,\beta):=\Card\{i:\alpha_i\neq \beta_i\}/2.
\]
\end{defn}

Recall that to a path $\beta$ can also be associated a link pattern, that is
a set of pairings between indices (possibly leaving some of them unpaired).
Each such pairing corresponds to a $\ldots,2,\ldots,1,\ldots$ in
the corresponding binary string.
Define a set of paths by $\mathcal{F}(\beta)$ 
as follows:
\[
\mathcal{F}(\beta):=\{\alpha\le\beta: \mbox{Some pairings of $\beta$ are flipped} \}.
\]
where by flipped we mean replacing $\ldots,2,\ldots,1,\ldots$ with $\ldots,1,\ldots,2,\ldots$ in the binary string of $\beta$.
If the number of pairings of $\beta$ is $r$, 
then the cardinality of $\mathcal{F}(\beta)$ is $2^r$. 

\begin{example}
$\beta=2121$, that
is the link pattern $\raisebox{0.1cm}{\linkpattern[inverted]{1/2,3/4}}$. 
$\mathcal{F}(\beta)=\{2121,1221,2112,1212\}$.
\end{example}

\begin{remark}\label{remark-f-II}
The set $\mathcal{F}(\beta)$ can be rephrased in terms of 
Dyck strips. Let us fix a path $\beta$. 
$\mathcal{F}(\beta)$ is the set of paths $\{\alpha: \alpha\le\beta\}$ (with
the $-$ convention) such that 
the region between them, denoted (following the notation of skew Ferrers
diagrams) $\beta/\alpha$, is filled with Dyck 
strips according to Rule II.
\end{remark}

Note that when $\alpha\in\mathcal{F}(\beta)$, $d(\alpha,\beta)$ is equal to the number 
of flipped pairings in $\beta$. 

On $\mathcal{M}^-_{N,K}$, let us define 
\[
\widetilde{C}_{\beta}
:=
\sum_{\alpha\in\mathcal{F}(\beta)}t^{-d(\alpha,\beta)}m_{\alpha}
\]
%where paths $\alpha,\beta\in\mathcal{P}_{N,K}$ correspond to 
%$x,y\in\mathcal{S}_N/(\mathcal{S}_K\times\mathcal{S}_{N-K})$.

Let $s_{i_1}s_{i_2}\ldots s_{i_l}=\short{\beta}$ be a reduced word of
the shortest coset representative $\short{\beta}$. 
We denote this ordered product by $\overleftarrow{\prod}_{s_i\in s(\beta)}s_i$.
\begin{prop}[see also \cite{KL-KL}]\label{dual-KL}
The basis $(\widetilde{C}_{\beta})$ for $\beta\in\mathcal{P}_{N,K}$ consists of 
elements that may be factorized as 
\[
%\label{expand-fac}
\widetilde{C}_{\beta}
=
\left(\prod_{s_i\in s(\beta)}^{\longleftarrow} (T_{i}+t^{-1})\right)m_{\beta_0}.
\]
where $\beta_0=(1\ldots12\ldots2)\in\mathcal{P}_{N,K}$.
 \end{prop}

\begin{proof}
We prove the proposition by induction on $\beta$. 
We have $\widetilde{C}_{\beta_0}=m_{\beta_0}$ and 
$\widetilde{C}_{s_K.\beta_0}=m_{s_K.\beta_0}+t^{-1}m_{\beta_0}
=(T_K+t^{-1})m_{\beta_0}$.

Let $\beta,\beta'\in\mathcal{P}_{N,K}$ and 
suppose the statement holds true for all $\beta'<\beta$. 
Now let $s(\beta)=s_is(\beta')$ with $|\beta|=|\beta'|+1$. 
This condition is equivalent to $(\beta_i,\beta_{i+1})=(\beta'_{i+1},\beta'_i)=(1,2)$. 

Note that the contribution of a pairing to $m_{\alpha'}$ for $\alpha'\in\mathcal{F}(\beta')$
is independent of each other. 
Therefore, it is enough to check the action of $T_{i}+t^{-1}$ on 
a partial path of $\alpha'$ involving $\alpha'_i$ and $\alpha'_{i+1}$. 
We have three cases.
\begin{enumerate}[(i)]
\item Suppose $\alpha'_i=1$ is unpaired and $(\alpha'_{i+1},\alpha'_{j})=(2,1)$ is a pairing.  
In this case, $\alpha''_{i}=1$ holds true for all $\alpha''\in\mathcal{F}(\beta)$.  
\[
(T_i+t^{-1})(m_{\ldots\underline{12}\ldots1\ldots}+t^{-1}m_{\ldots\underline{11}\ldots2\ldots})
=m_{\ldots\underline{21}\ldots1\ldots}+t^{-1}m_{\ldots\underline{12}\ldots1\ldots}
\]
where $T_i$ acts on the underlined places. 
Now $\alpha_{j}=1$ becomes an unpaired $1$, and $(\alpha_i,\alpha_{i+1})$ 
becomes a pairing in $\alpha$.
Suppose $\alpha'_{i+1}=2$ is unpaired and $(\alpha'_j,\alpha_i)=(2,1)$ is a pairing. 
Similarly as above, we have $\alpha_j=2$ is unpaired and $(\alpha_i,\alpha_{i+1})=(2,1)$  
is a pairing.
\item $(\alpha'_k,\alpha'_i)=(\alpha'_{i+1},\alpha'_{l})=(2,1)$ with $k<i, i+1<l$ and they are pairings. 
\begin{align*}
(T_i+t^{-1})(m_{\ldots2\ldots\underline{12}\ldots1}+t^{-1}m_{\ldots1\ldots\underline{22}\ldots1}
+t^{-1}m_{\ldots2\ldots\underline{11}\ldots2}+t^{-2}m_{\ldots1\ldots\underline{21}\ldots2}) \\
=m_{\ldots2\ldots\underline{21}\ldots1}+t^{-1}m_{\ldots2\ldots\underline{12}\ldots1}
+t^{-1}m_{\ldots1\ldots\underline{21}\ldots2}+t^{-2}m_{\ldots1\ldots\underline{12}\ldots2},
\end{align*}
This implies that $(\alpha_{k},\alpha_{l})=(\alpha_{i},\alpha_{i+1})=(2,1)$ are pairings 
in $\alpha$. 
\item Suppose both $\alpha'_i=2$ and $\alpha'_{i+1}=1$ are unpaired. We have 
\begin{align*}
(T_{i}+t^{-1})m_{\ldots\underline{12}\ldots}=m_{\ldots21\ldots}+t^{-1}m_{\ldots21\ldots}, 
\end{align*}
which means $(\alpha_i,\alpha_{i+1})=(2,1)$ is a pairing. 
\end{enumerate}
In all cases, obtained expression gives us the set $\mathcal{F}(\beta)$ and 
desired coefficients.  
\end{proof}

\begin{prop}\label{theorem-minus}
The basis $(\widetilde{C}_\beta)$ is the Kazhdan--Lusztig basis $(C^-_{\beta})$. 
\end{prop}
\begin{proof}
Note that $\widetilde{C}_\beta$ is invariant under the bar involution
since $\overline{T_i+t^{-1}}=T_i+t^{-1}$ and 
$\overline{m_{\alpha_0}}=m_{\alpha_0}$. 

From Proposition~\ref{dual-KL}, it is clear that the coefficient 
for $m_\beta$ is $1$ and all other coefficients 
for $m_{\alpha}$ are in $t^{-1}\mathbb{N}[t^{-1}]$ for $\alpha<\beta$.
\end{proof}
When the region $\beta/\alpha$ is filled with 
Dyck strips via Rule II, it is clear that $d(\alpha,\beta)$ is equal to 
the number of Dyck strips. 
From Proposition~\ref{theorem-minus} together with Remark~\ref{remark-f-II}, 
we have the (see also \cite{Bre02})
\begin{theorem}\label{cor-minus-KL}
The Kazhdan--Lusztig polynomial $P^-_{\alpha,\beta}$ is given by 
\[
 P^-_{\alpha,\beta}(t^{-1})=Q^{II,-}_{\alpha,\beta}(t^{-1})=
t^{-d(\alpha,\beta)}, \qquad \alpha\in\mathcal{F}(\beta).
\]
\end{theorem}

As mentioned in the introduction, the parabolic KL basis of
$\M^-$ is closely related to the formulation of
the Temperley--Lieb algebra in terms of tangles as used
in knot theory \cite{Kauf-TL}. Indeed, in this basis,
the operators $T_i+t^{-1}$ which appeared in the proof,
and which are the usual generators in terms of which the Temperley--Lieb
algebra is formulated, have a natural graphical action on link patterns:
they correspond to pasting a 
$\vcenter{\hbox{\begin{tikzpicture} 
\tikzset{vertex/.style={circle,draw=black,line width=0.3pt,fill=red,inner sep=1.5pt}}
\tikzset{edge/.style={bend left=60,draw,very thick,blue}}
\draw[edge] (0,0) coordinate[vertex] to (1,0) coordinate[vertex];
\draw[edge] (1,0.7) coordinate[vertex] to (0,0.7) coordinate[vertex];
\end{tikzpicture}}}$
to the link pattern, i.e.,
reconnecting the existing pairings between
neighboring sites $i$ and $i+1$ and creating a new pairing $(i,i+1)$.

Descriptions of $P^-$ that are analogous to Thm.~\ref{cor-minus-KL}
appear under various guises in the literature; see \cite[Eq.~(5.12)]{BS3} 
for an alternative form of it in terms of oriented cup diagrams,
\cite[Lemma 2.2]{SW-springer} for an interpretation of this formula
in terms of Springer fibres; % \cite[Thm.~5.4]{Bs3}, \cite[Thm.~5.2]{Sperv} 
%for a categorical interpretation of the KL-polynomials in terms of graded
%decomposition numbers in highest weight categories;
\cite[Sect.~8]{MNdGB} for an appearance in statistical loop models;
and \cite[Sect.~8]{Martin-Brauer} for a connection to the representation
theory of the Brauer algebra.

\subsection{The inversion formula}
In preparation for the study of the module $\mathcal{M}^+_{N,K}$,
we invert the matrix $Q^{II,-}$.

\begin{theorem}\label{mainthm}
\[
\sum_{\beta\in\mathcal{P}_{N,K}} Q_{\alpha,\beta}^{I,-}(t^{-1}) Q^{II,-}_{\beta,\gamma}(t^{-1}) (-1)^{|\beta|+|\gamma|}
=\delta_{\alpha,\gamma}
\]
\end{theorem}
\begin{proof}
If $\alpha\not\le\gamma$ the l.h.s.\ is zero, and if $\alpha=\gamma$ 
it is one. We now assume $\alpha<\gamma$.
By definition,
\newcommand\DS{\mathcal{D}}
\newcommand\DSI{\mathcal{D}^{I}}
\newcommand\DSII{\mathcal{D}^{II}}
\[
\sum_\beta Q_{\alpha,\beta}^{I,-}(t^{-1}) Q^{II,-}_{\beta,\gamma}(t^{-1}) (-1)^{|\beta|+|\gamma|}
=\sum_{\beta,\alpha\le\beta\le\gamma}
\sum_{\DSI\in\CI(\alpha,\beta)}
\sum_{\DSII\in\CII(\beta,\gamma)}
t^{-(|\DSI|+|\DSII|)} (-1)^{|\DSII|}
\]
The sign was obtained by noting
that Dyck strips have odd length, so that the number of boxes $|\gamma|-|\beta|$
and the number of Dyck strips $|\DSII|$ of $\DSII$ have same parity.
Now merge together the two families of Dyck strips $\DSI$
and $\DSII$ into a single family $\DS$, and switch the summations:
\[
\sum_\beta Q_{\alpha,\beta}^{I,-} Q^{II,-}_{\beta,\gamma} (-1)^{|\beta|+|\gamma|}
=\sum_{\DS\in{\rm Conf}(\alpha,\beta)}
t^{-|\DS|} 
\sum_{\beta\in P(\DS)}
(-1)^{|\DSII(\beta)|}
\]
where $P(\DS)$ is the set of paths $\beta$ between $\alpha$ and $\gamma$
such that the $D\in\DS$ below $\beta$ satisfy rule I and
those above $\beta$ satisfy rule II; we denote the corresponding partition
$\DS=\DSI(\beta)\sqcup \DSII(\beta)$. 

We shall show that for a fixed decomposition $\DS$ of 
$\gamma/\alpha$ into Dyck strips, the sum over $\beta$, i.e.,
over subdivisions of $\DS$ into two subsets 
(one satisfying rule I, the other rule II), is zero.
In all that follows, we assume $P(\DS)\ne\emptyset$ (otherwise
the sum is trivially zero). 

In this proof we shall need a relation on Dyck strips
in $\DS$, which mimics the definition of rule II.
We write 
that $D\prec D'$ if all boxes just above, NW or NE of a box of $D$ belong
to $D$ or $D'$. This relation has a tree structure in the sense
that for given $D$ there is at most one $D'$ such that $D\prec D'$.
If there are no such $D'$, then $D$ is called a minimal element
(this is just the usual notion of minimality for the associated
order relation).

Define
\[
\mathcal{I}(\DS)=
%\{ D: \exists \beta,\beta'\in P(\DS): D\in \DSI(\beta)\cap \DSII(\beta')\}
\Big(\bigcup_{\beta\in P(\DS)} \DSI(\beta)\Big)
\cap
\Big(\bigcup_{\beta\in P(\DS)} \DSII(\beta)\Big)
\]
i.e., the set of
Dyck strips which can be on either side of the boundary between
zones I and II.
We have the first observation
\begin{lemma}
If $D,D'\in \mathcal{I}(\DS)$, $D\ne D'$,  
then the $x$ coordinates of boxes
of $D$ and $D'$ are distant by at least 2. 
\end{lemma}
\begin{proof}
Assume the $x$ coordinates are distant by less than 2.
Then there is a box $(x,y)$ of one of the two Dyck strips, 
say $D$, which is above a box $(x',y')$ of $D'$ in the sense
that $y>y'$ and $x=x'\pm 1$. But note that this excludes the possibility
of finding a path between $\alpha$ and $\gamma$ such that
$D$ is below it and $D'$ is above it.
Therefore, choosing $\beta,\beta',\beta'',\beta'''$ such that 
$D\in \DSI(\beta)\cap \DSII(\beta')$ and $D'\in \DSI(\beta'')\cap \DSII(\beta''')$,
we conclude that $D'\in \DSI(\beta)$ and $D\in \DSII(\beta''')$.
But this implies that there is a region containing
both $D$ and $D'$, namely the region below $\beta$ and above $\beta'''$, 
in which both rule I and rule II apply.
The rule II and the relative position of $D$ and $D'$
imply that there is a chain $D'\prec D_1\prec\cdots \prec D_k\prec D$;
but this in turn implies that 
two successive Dyck strips in the chain also satisfy the conditions of
applicability of rule I.
These two facts are contradictory because they imply
opposite inequalities on the lengths, i.e., 
$l(D')<l(D_1)<\cdots<l(D_k)<l(D)$
and $l(D')\ge l(D_1)\ge\cdots\ge l(D_k)\ge l(D)$.
\end{proof}

We conclude immediately that distinct elements of $\mathcal{I}$ ``do not interact'' with
each other in the sense that they can be added/removed independently from $\DS^{I}$, $\DS^{II}$.
More explicitly, note that since $P(\DS)\ne\emptyset$,
$\bigcup_{\beta\in P(\DS)} \DSII(\beta)\ne\emptyset$; and its
lower boundary is again a path, say $\alpha_0$. 
Similarly one can define $\gamma_0$ which is the upper boundary of 
$\bigcup_{\beta\in P(\DS)} \DSI(\beta)$.
Then for any subset $\mathcal{J}\subset \mathcal{I}(\DS)$,
there is a path $\beta\in P(\DS)$ such that
$\DSII(\beta)=\DSII(\gamma_0)\sqcup \mathcal{J}$ and
$\DSI(\beta)=\DSI(\alpha_0)\sqcup (\mathcal{I}(\DS)\backslash\mathcal{J})$.
Indeed, rules I and II cannot apply to two elements of $\mathcal{I}(\DS)$ because they
are too far apart, and in all other cases one easily checks that these rules are already
satisfied by definition.

To summarize, we have found that the summation over $\beta$
is structured as follows: $P(\DS)$ is of cardinality $2^{|\mathcal{I}(\DS)|}$,
corresponding to whether $D\in\mathcal{I}(\DS)$, 
is above or below the path separating regions I and II.
Furthermore, we have the following key fact:
\begin{prop}\label{keyprop}
\[
\mathcal{I}(\DS)\ne\emptyset
\]
\end{prop}
\begin{proof}
We shall in fact provide an explicit description of $\mathcal{I}(\DS)$
using the relation $\prec$.
Recall from the structure of $P(\DS)$ described above
that there is a path $\alpha_0\in P(\DS)$ such that
$\bigcup_{\beta\in P(\DS)} \DSII(\beta)=\DSII(\alpha_0)$.

We claim that $\mathcal{I}(\DS)$ is exactly the set 
$\min \DSII(\alpha_0)$ of minimal elements (in the sense of $\prec$) 
of $\DSII(\alpha_0)$. $\mathcal{I}(\DS)\subset \min \DSII(\alpha_0)$ 
is evident by definition of $\mathcal{I}(\DS)$. %$\DSII(\alpha_0)$.
Let us now prove the reverse inclusion, i.e., prove that any minimal element
of $\DSII(\alpha_0)$ can also be moved to the region I.

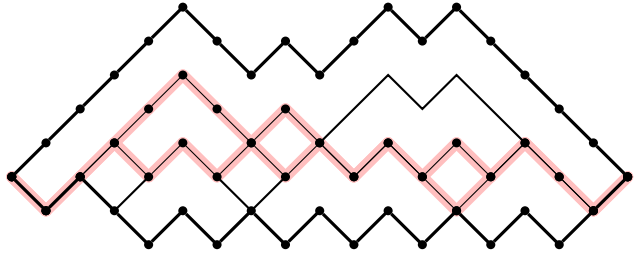
\begin{figure}[h!]
\setlength{\unitsize}{0.45cm}
\begin{tikzpicture}[x=\unitsize,y=\unitsize,baseline=0]
\tikzset{vertex/.style={}}%
\tikzset{edge/.style={pink,line width=0.15cm,line cap=round}}%
\rawpath{-,+,+,-,+,-,+,-,+,-,+,-,-,+,+,-,-,+}
\rawpath{-,+,+,+,+,-,-,+,-,-,+,-,+,-,+,-,-,+}
\tikzset{vertex/.style={circle,fill=black,inner sep=1.2pt,outer sep=0pt}}%
\tikzset{edge/.style={very thick}}%
\rawpath{+,+,+,+,+,-,-,+,-,+,+,-,+,-,-,-,-,-}
\rawpath{-,+,-,-,+,-,+,-,+,-,+,-,+,-,+,-,+,+}
\tikzset{edge/.style={}}
\rawpath{-,+,+,-,+,-,+,-,+,-,+,-,-,+,+,-,-,+}
\rawpath{-,+,+,+,+,-,-,+,-,-,+,-,+,-,+,-,-,+}
\draw[thick] (3,-1) -- (4,0) (6,0) -- (7,-1) -- (8,0) 
(9,1) -- (11,3) -- (12,2) -- (13,3) -- (15,1);
\end{tikzpicture}%
\caption{Illustration of the proof of Prop.~\ref{keyprop}. The thick lines represent
the boundaries $\alpha_0$ and $\gamma_0$ of the maximal/minimal regions I/II, 
so that the Dyck strips in between the two form $\mathcal{I}(\DS)$.}
\end{figure}
\rem{also maybe a figure explaining the contradiction in the proof?}

Pick such a minimal element $D_{min}\in \DSII(\alpha_0)$. 
Due to the way we defined $\prec$, it is easy to see that
$\DSII(\alpha_0)\backslash D_{min} = \DSII(\beta)$ 
for some $\beta$ above $\alpha_0$.
We now claim that the Dyck strips below $\beta$ satisfy rule I.
These Dyck strips consist of the Dyck strips below $\alpha_0$,
which by construction satisfy rule I, plus the additional $D_{min}$.
To a box of $D_{min}$ with coordinates $(x,y)$ associate $D_x$, which is 
the Dyck strip to which belongs the box right below, i.e., $(x,y-2)$,
or $\varnothing$ if this box is below the bottom line $\alpha$. 
Rule I means that this function should be constant. The
proof is by contradiction. Suppose
there is an $x$ such that $D_x\ne D_{x+1}$. We can assume up to reflection
w.r.t.\ the $y$ axis that the higher of the two boxes is the first,
i.e., $(x,y-2)\in D_x\ne\varnothing$ and $(x+1,y-3)\in D_{x+1}$ if
$D_{x+1}\ne\varnothing$.
Now note that
$y-3\ge h-2$, where $h$ is the minimum $y$ coordinate of boxes of $D_{min}$; 
so that $y-2\ge h-1$. 
Therefore the Dyck strip $D_x$ cannot pass below
the endpoints of $D_{min}$ (whose $y$ coordinates are $h$); in other words,
its $x$ span is strictly included in that of $D_{min}$ and it touches
$D_{min}$ at its two boundaries. 
%Let us call
%this span $\{x_1,x_1+1,\ldots,x_2\}$; the existence of $D_x\ne\varnothing$
%implies that $D_{x'}\ne\varnothing$
%for all $x_1\le x'\le x_2$.

%To conclude, note that the preimages of the 
%map $x\mapsto D_x$ provide a set partition
%of $\{x_1,\ldots,x_2\}$ which is {\em non-crossing}\ since Dyck strips
%cannot cross each other. Every non-crossing partition possesses
%a block that is reduced to a single interval (easy proof by induction
%on the width). By inspection, one checks that the 
%corresponding Dyck strip $D_{x'}$ is such that $D_{x'}\prec D$,
%contradicting the minimality of $D$.

Now introduce another relation $\to$ as follows:
$D\to D'$ if {\em there exists}\/ a box of $D'$ NW, NE or above a box of $D$.
There is a naturally associated order relation, which we simply denote 
$D\to\cdots\to D'$, obtained by forming chains.
We can consider
$\mathcal{X}=\{ D\in\DS^I(\alpha_0): D_x\to\cdots\to D \}$.
$\mathcal{X}$ is non-empty because $D_x\in \mathcal{X}$; 
a maximal element $D$ of it
(for the order relation $\to\cdots\to$)
is such that all boxes NW, NE and above it are outside $\DS^I(\alpha_0)$;
but since its $x$ span is strictly included in that of $D_{min}$, these boxes
must belong to $D_{min}$.
Therefore $D\prec D_{min}$, contradicting the minimality of $D_{min}$.

Thus, since
$\DSII(\alpha_0)\ne\emptyset$, $\mathcal{I}(\DS)=\min \DSII(\alpha_0)\ne\emptyset$.
\end{proof}

Note finally that the possible paths $\beta\in P(\DS)$ correspond
to $\DS^{II}(\beta)=\DS^{II}(\gamma_0)\sqcup \mathcal{J}$ where 
$\mathcal{J}$ is any subset of $\mathcal{I}(\DS)$, so that
we can compute the sum over $\beta\in P(\DS)$ by rewriting it 
$(-1)^{|\DS^{II}(\gamma_0)|}\sum_{\mathcal{J}\subset\mathcal{I}(\DS)}(-1)^{|\mathcal{J}|}=0$.
\end{proof}

\rem{corollary for minus case? probably don't really care}

\subsection{On the module $\mathcal{M}^+_{N,K}$}
The two families 
of Kazhdan--Lusztig polynomials $P^\pm_{x,y}(t^{-1})$ on the modules 
$\mathcal{M}^\pm_{N,K}$ are related by 
the duality theorem~\ref{paraduality}.
Together with Lemma~\ref{lemma-flip}, we have 
\begin{align*}
(P^+)^{-1}_{\alpha,\beta}(t^{-1})&=(-1)^{|\alpha|+|\beta|}P^-_{\beta,\alpha}(t^{-1}) \\
&=Q^{II,+}_{\alpha,\beta}(-t^{-1}),
\end{align*}
where we have once again used that $(-1)^{|\alpha|+|\beta|}=(-1)^{|\alpha|-|\beta|}$ and 
that the length of a Dyck strip is always odd.
% Taking into account the change of convention $\epsilon$, we define:
% \begin{defn}
% A set of paths $\mathcal{L}(\beta)$ is defined as 
% %
% \begin{align}
% \mathcal{L}(\beta):=\{\alpha : \beta\in\mathcal{F}'(\alpha) \}  
% \end{align}
% %
% where $\mathcal{F}'(\alpha)$ is the same as $\mathcal{F}(\alpha)$ 
% except flipping $(1,2)$ instead of $(2,1)$. 
% \end{defn}
%
%
%
% \begin{example}
% Let $\beta=2121$. Then $\mathcal{L}(\beta)=\{2121, 1221, 2112, 1212,  1122\}$.
% \end{example}
%
% \begin{remark}\label{remark-L-II}
% The set $\mathcal{L}(\beta)$ can be rephrased in terms of the Ferrers diagram and 
% Dyck strips. Let us fix a path $\beta$. 
% Take the (45 degrees rotated) Ferrers diagram $\lambda(\beta)$ with 
% $\epsilon=+$ convention.  
% $\mathcal{L}(\beta)$ is a set of paths $\{\alpha: \alpha\le\beta\}$ such that 
% a skew Ferrers diagram $\lambda(\beta)/\lambda(\alpha)$ is filled with Dyck 
% strips by Rule II.
% \end{remark}
% %
% Recall that Lemma~\ref{lemma-flip} deals with the change of conventions $\epsilon$ in 
% the pictorial expressions.  
% However, when we consider paths as binary strings, we have to be careful of this change 
% since $1$ and $2$ are exchanged.
% This is why small complication occurs in the definition of $\mathcal{L}(\beta)$
% compared to $\mathcal{F}(\beta)$.
%
%
% Since we have the ``transposed'' set of $\mathcal{F}(\beta)$, we have immediately
%
Hence,
\begin{cor}\label{inverse-LS}
On $\mathcal{M}^+_{N,K}$, the monomial basis $m_\alpha$ is expressed 
in terms of the Kazhdan--Lusztig basis as 
\[
 m_\beta
 =
\sum_{\alpha\le\beta}Q^{II,+}_{\alpha,\beta}(-t^{-1})C^+_{\alpha}.
\]
\end{cor}
A slightly more explicit version of this formula
is provided in appendix C.

More importantly, Theorem~\ref{mainthm} allows us to invert this relation
to obtain the Kazhdan--Lusztig basis $C^+_\beta$ in terms of 
the monomial basis $m_\alpha$:
\begin{cor}\label{theorem-LS}
The Kazhdan--Lusztig polynomial $P^+_{\alpha,\beta}$ 
is the generating function of Dyck strips according to Rule I, that is,
\[
 P^+_{\alpha,\beta}(t^{-1})=Q^{I,+}_{\alpha,\beta}(t^{-1}).
\]
\end{cor}
Examples can be found in appendix A.

Formulae for such polynomials were of course already known:
see in particular \cite{ES-KL} for a similar approach in a more
general setting; and \cite{LS-KL}, the combinatorial description of which
is described in the next section and shown to be in bijection with ours.

\section{Relation to the Lascoux--Sch\"utzenberger rule}

\subsection{Lascoux--Sch\"utzenberger binary trees}
We briefly review the construction of the binary trees of 
Lascoux--Sch\"utzenberger 
to compute the Kazhdan--Lusztig polynomials for Grassmannian permutations \cite{LS-KL}. 
In our setup, they correspond to the polynomials $P^+_{\alpha,\beta}(t^{-1})$ 
from Corollary~\ref{theorem-LS}.  

Let $\mathcal{Z}$ be a set such that 
$\emptyset\in\mathcal{Z}$ (where $\emptyset$ represents the empty string), 
$z\in\mathcal{Z}\Rightarrow1z2\in\mathcal{Z}$ and 
$z_1,z_2\in\mathcal{Z}\Rightarrow z_1z_2\in\mathcal{Z}$. 
We define inductively 
a rooted tree $A(w)$ for $w$ an arbitrary binary string by:
\begin{itemize}
\item $A(\emptyset)$ is the empty tree,
 \item $A(2w)=A(w1)=A(w)$,
 \item $A(zw)$ is obtained by attaching the trees for $A(z)$ and $A(w)$ at their roots, $z\in\mathcal{Z}$,
 \item $A(1z2)$ is obtained by attaching an edge just below the tree $A(z)$,
$z\in\mathcal{Z}$.
\end{itemize}

We denote by $\|\alpha\|$ the length of a binary string $\alpha$ and 
by $\|\alpha\|_{\sigma}$ the number of $\sigma$ in a string $\alpha$.
Let $v,w\in\{1,2\}^{N}$ with $v\le w$. 
and $v=v'\alpha\beta v'', w=w'\underline{12}w''$ with $\|\alpha'\|=\|w'\|$ and $\alpha,\beta\in\{1,2\}$. 
A {\em capacity}\/ of the edge corresponding to the underlined $1$ and $2$ is defined by
\begin{align}
 \mathrm{cap}(12):=\|v'\alpha\|_1-\|w'1\|_1
\end{align}
The condition $v\le w$ implies a capacity is always non-negative.

The capacity of $v$ with respect to $w$ is the collection of capacities of pairs of 
adjacent $1$ and $2$ in $w$ and called the relative capacities. 

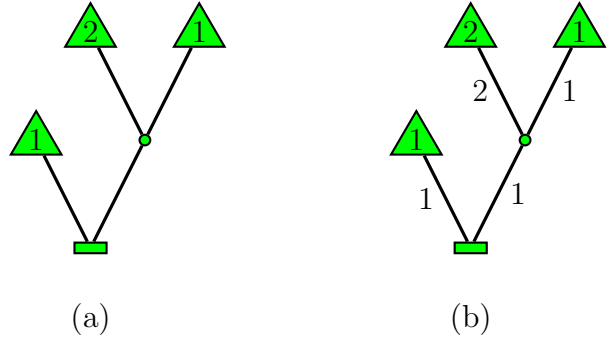
\begin{figure}[ht!]
\begin{tikzpicture}[x=\unitsize,y=\unitsize,label distance=5mm]
%\tikzset{vertex/.style={circle,draw=black,thick,fill=green,inner sep=0pt,minimum size=4pt}}
\treestyle
\node[root,label=below:(a)] {} [tree]
child {node[leaf] {1}}
child {node[vertex] {}
  child {node[leaf] {2}}
  child {node[leaf] {1}}
};
\begin{scope}[xshift=5cm]
\node[root,label=below:(b)] {} [tree]
child {node[leaf] {1} 
       edge from parent node[left]{1}}
child {node[vertex] {} 
  child {node[leaf] {2} 
         edge from parent node[left]{2}}
  child {node[leaf] {1} 
         edge from parent node[right]{1}}
  edge from parent node[right]{1}
};
\end{scope}
\end{tikzpicture}%
\caption{(a) A tree with capacities. (b) A labelled tree.}\label{figtree}
\end{figure}

We denote by $A(w/v)$ the rooted tree with relative capacities.
$A(w/v)$ is obtained from the tree $A(w)$ by putting corresponding
capacities at leaves (end points) of the tree, see Fig.~\ref{figtree}(a). 

A {\em labelling}\/ of the tree $A(w/v)$ is a set of non-negative 
integers on edges of $A(w)$ satisfying
\begin{itemize}
 \item An integer on a leaf is less than or equal to its capacity, 
 \item Integers on edges are non-increasing from leaves to the root.
\end{itemize}
See Fig.~\ref{figtree}(b).

The analysis of the recursive relations for both the Kazhdan--Lusztig 
polynomials and the generating function of the tree $A(w/v)$ led
Lascoux and Sch\"utzenberger to the following theorem, 
formulated here in our notations
(in particular we identify as before binary strings and paths 
with convention $+$):
\begin{theorem}[Lascoux, Sch\"utzenberger]\label{lsthm}
\[
 P^{+}_{\alpha,\beta}(t^{-1})=t^{|\alpha|-|\beta|}\sum_{\nu}t^{2\Sigma(\nu)}.
\]
where $\nu$ runs over all possible labellings of $A(\beta/\alpha)$ 
and $\Sigma(\nu)$ is the sum of labels of $\nu$.
\end{theorem}

Below, we produce a bijection between a labelling of $A(\beta/\alpha)$ and 
a configuration of Dyck strips between paths $\alpha$ and $\beta$
(i.e., in the skew-diagram $\beta/\alpha$).

\subsection{From trees to link patterns}
\begin{figure}[h!]
\setlength{\unitsize}{0.75cm}%
{
\begin{tikzpicture}[x=\unitsize,y=\unitsize,baseline=0]
%\tikzset{vertex/.style={circle,draw=black,fill=red,inner sep=1.5pt}}%
%\tikzset{edge/.style={bend right=75,draw,very thick,blue,looseness=1.5}}%
\tikzset{edge/.append style={bend left=75,min distance=0.8cm,looseness=1.2}}%
%\numberedfalse%
\linkpattern[unit=\unitsize,tikzstarted,inverted]{1/2,3/8,4/5,6/7}%
\treestyle
\node[root] (a) at (3.5,-2.7) {};
\node[leaf] (b) at (1.5,-0.35) {};
\node[vertex] (c) at (5,-1.3) {};
\node[leaf] (d) at (4.5,-0.35) {};
\node[leaf] (e) at (6.5,-0.35) {};
\draw[very thick] (a) -- (b.south) (a) -- (c) -- (d.south) (c) -- (e.south);
\end{tikzpicture}%
}
\qquad%
\begin{tikzpicture}[x=\unitsize,y=\unitsize,baseline=0]
%\tikzset{vertex/.style={circle,draw=black,fill=red,inner sep=1.5pt}}%
%\tikzset{edge/.style={bend right=75,draw,very thick,blue}}%
\numberedfalse%
\linkpattern[unit=\unitsize,tikzstarted,inverted]{1/2,3/8,4/5,6/7}%
\node at (1.5,-0.75) {1};
\node at (4.5,-0.75) {1};
\node at (6.5,-0.75) {0};
\node at (5.5,-1.8) {1};
\end{tikzpicture}%
\\[4pt]
(a)\hskip5cm(b)
\caption{(a) Link pattern and binary tree. (b) Labelling of the link pattern.}\label{linkpatterntotree}
\end{figure}
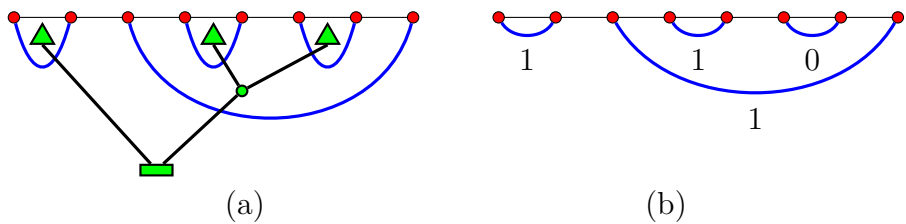

In the previous section we have introduced, 
following Lascoux and Sch\"utzenberger,
binary trees starting from a binary string. Using
the bijections of section \ref{secbij}, we can equivalently start from a path,
or from a link pattern. The latter correspondence is particularly natural,
since the binary tree is the {\em dual graph}
of the link pattern, cf Fig.~\ref{linkpatterntotree}(a) (with the same example as in 
Fig.~\ref{figtree}).
Note that there is a bijective map $p$ which to an edge $e$ associates a pairing
$p(e)$ of the link pattern.
However, unless the link pattern has maximum number of pairings,
the map from link patterns to trees is not one-to-one:
when we take the dual graph, we ignore the unpaired vertices.
In what follows we denote by $\pi(\beta)$ the link pattern associated to
the path (or binary string) $\beta$.

%We have a construction of a tree $A(\beta)$ from a path 
%(or equivalently binary code) $\beta$.
%Inverting this construction, we get a path $\beta'$ for $A(\beta)$, 
%which can be different from $\beta$.  
%Also we have a bijection from $\beta$ and $\beta'$ to link patterns 
%$\pi(\beta)$ and $\pi(\beta')$.
%However, these two link patterns contain the same data for planar pairings 
%in the sense that we get $\pi(\beta')$ by removing the unpaired vertices 
%from $\pi(\beta)$.
%Therefore, hereafter, we identify both $\pi(\beta)$ and $\pi(\beta')$ as a 
%link pattern.  
%Each planar pairing of the link patten $\pi(\beta)$ is exactly corresponding 
%to a unique edge of the tree $A(\beta)$.  
%A link pattern is obtained by taking ``dual" picture of a tree.
%[See Fig~..]\\
%

It is also convenient to attach labellings to the link pattern as follows.
Given a labelling of $A(\beta/\alpha)$ and an edge $e$ with label $n(e)$,
we put the label $n'(p(e))=n(e)-n(e')$ on the corresponding pairing $p(e)$, where 
the edge $e'$ is the parent edge of $e$, unless there is no parent edge (edge
connected to the root) in which case we put $n(e)$. See 
Fig.~\ref{linkpatterntotree}(b) (with the same labelling as in 
Fig.~\ref{figtree}(b)).

Labellings of the link pattern $\pi(\beta)$ thus obtained
from a labelling of $A(\beta/\alpha)$ are defined by the following conditions:
\begin{itemize}
 \item All labels are non-negative integers.
 \item Given a smallest planar pairing $p(e)$ (a pairing of neighbors), 
the sum of all labels on planar pairings 
which surround $p(e)$ is less than or equal to the capacity of $e$.  
\end{itemize}

\subsection{From labelled link patterns to Dyck strips}
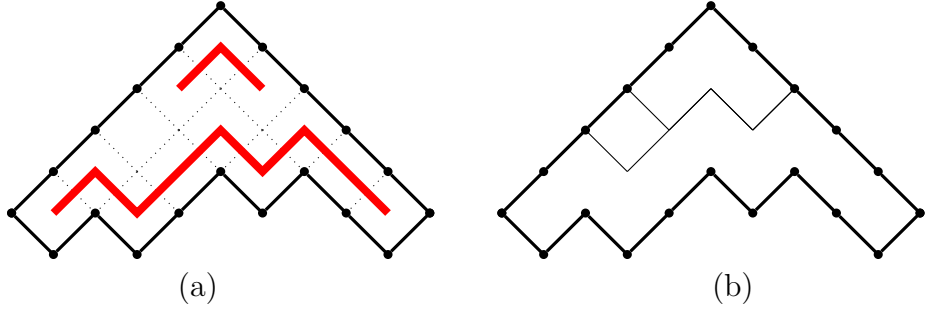
\begin{figure}[h!]
\setlength{\unitsize}{0.55cm}
\begin{tikzpicture}[x=\unitsize,y=\unitsize,baseline=0]
\tikzset{vertex/.style={circle,fill=black,inner sep=1.2pt,outer sep=0pt}}%
\tikzset{edge/.style={very thick}}%
\rawpath{+,+,+,+,+,-,-,-,-,-}%
\rawpath{-,+,-,+,+,-,+,-,-,+}%
\tikzset{vertex/.style={}}%
\tikzset{edge/.style={dotted}}%
\rawpath{+,+,+,+,-,+,-,-,-,-}%
\rawpath{+,+,+,-,+,-,+,-,-,-}%
\rawpath{+,+,-,+,-,+,-,+,-,-}%
\rawpath{+,-,+,-,+,-,+,-,+,-}%
\draw[red, line width=3pt] (1,0) -- (2,1) -- (3,0) -- (5,2) -- (6,1) -- (7,2) -- (9,0) (4,3) -- (5,4) -- (6,3);
\end{tikzpicture}%
\qquad
\begin{tikzpicture}[x=\unitsize,y=\unitsize,baseline=0]
\tikzset{vertex/.style={circle,fill=black,inner sep=1.2pt,outer sep=0pt}}%
\tikzset{edge/.style={very thick}}%
\rawpath{+,+,+,+,+,-,-,-,-,-}%
\rawpath{-,+,-,+,+,-,+,-,-,+}%
\tikzset{vertex/.style={}}%
\tikzset{edge/.style={}}%
\rawpath{+,+,+,-,+,-,+,-,-,-}%
\rawpath{+,+,-,+,+,-,+,-,-,-}%
\end{tikzpicture}%
\\[4pt]
(a)\hskip6.5cm(b)
\caption{Dyck paths and Dyck strips.}\label{lptods}
\end{figure}

We now consider a pair of paths $\alpha$ and $\beta$,
with $\alpha$ above $\beta$, and the associated link pattern $\pi(\beta)$ 
along with a labelling as above.
We associate to it a collection of Dyck strips between paths $\alpha$ and $\beta$
as follows.
Recall that a Dyck strip is characterized by a Dyck path. To each pairing $p$ of $\pi(\beta)$
we associate Dyck paths which start a half-step to the left of the left point of the pairing
and a half-step to the right of its right point. More precisely, if $p$ has label $n'(p)$ we then
stack $n'(p)$ such Dyck paths on top of each other, forming parallel layers above $\beta$. 
We then repeat the process
for every pairing, respecting the order which is to start with the largest arches
and end with the smallest arches (this way we respect rule I). 
See Fig.~\ref{lptods}(a) for the same example as in previous figures.
Note that some Dyck paths may have coinciding starting or end points, in which case
they are merged into a larger Dyck path. 

To each Dyck path (where Dyck paths which touch have been merged) we now associate
the corresponding Dyck strip. Note that such strips necessarily have length greater or
equal to 3. We claim that these strips remain under the path $\alpha$.
Indeed, let $p$ be a smallest planar pairing, that is, connecting $i$ and $i+1$. 
%Let $c'(p)$ be the sum of all non-negative integers on planar pairings, which 
%contain the pairing $p$ inside of them. 
Then the difference of heights of $\alpha$ and $\beta$ at the center of the pairing
(i.e., the depth of the $\wedge$-corner in the skew Ferrers diagram) is by direct
computation exactly the capacity of the corresponding edge $e$ in the tree $A(\beta/\alpha)$. 
Therefore the number of Dyck strips above that point, that is the sum of labels of
pairings surround $p$, which is nothing but the label of $e$ in the tree $A(\beta/\alpha)$,
is less or equal to the capacity, i.e., the difference of heights. Therefore the Dyck strips
remain below $\alpha$ at every local maximum of $\beta$, therefore everywhere.
%Therefore, it is clear that 
%\begin{align}
%\max\{c'(p) : \mathcal{D}\in\mathrm{Conf}(\alpha,\beta)\}=\mathrm{Cap}(e).
%\end{align}
 
The last stage is to declare that the boxes of $\beta/\alpha$ that
do not belong to any of the Dyck strips above are by definition Dyck
strips of length one (consisting of a single box). See Fig.~\ref{lptods}(b) for the final result.

It is easy to show that the correspondence above is bijective. Therefore we have
proved the
\begin{theorem}
There exists a bijection between labellings of the tree $A(\beta/\alpha)$ and 
configurations of Dyck strips between paths $\alpha$ and $\beta$ satisfying rule I.
\end{theorem}

In order to show that Corollary \ref{theorem-LS} and Theorem
\ref{lsthm} are equivalent, we still
need to compare powers of $t^{-1}$, which naively look quite different. Let us start from a
configuration of Dyck strips between paths $\alpha$ and $\beta$.
Consider a Dyck strip of length $\ge3$. It 
is obtained from one or possibly several Dyck paths each associated to a certain
pairing, say $p_1,\ldots,p_k$. The number of boxes of the Dyck strip is equal to
$1+\sum_{i=1}^k d(p_i)$, where $d(p)$ is the distance between the
two endpoints of $p$.
This formula still holds for Dyck strips of size one provided
one associates to it zero Dyck paths. We now write the number of boxes between $\alpha$
and $\beta$ as:
\begin{align*}
|\beta|-|\alpha|&=\sum_{\substack{\text{Dyck strip}\\\text{formed from paths}\\p_1,\ldots,p_k}} 
\Big(1+\sum_{i=1}^k d(p_i)\Big)\\
&=\text{number of Dyck strips}+\sum_{\text{$p$ pairing}} n'(p) d(p)\\
&=\text{number of Dyck strips}+\sum_{\text{$e$ edge}} (n(e)-n(e')) (2+2\times\text{number of descendents of $e$)}\\
&=\text{number of Dyck strips}+2\sum_{\text{$e$ edge}} n(e)
\end{align*}
where we have used the fact that $d(p)=2+2$ 
times the number of pairings surrounded by $p$ and
translated it into the language of trees. We then write $n'(p(e))=n(e)-n(e')$ 
where $e'$ is the parent of $e$.
The final equality provides the required
identification of powers of $t^{-1}$.

Appendix B provides the full computation of a KL polynomial in the
various formulations (path, tree).

\appendix

\section{Table of polynomials at $N=4$, $K=2$}
\setlength{\unitsize}{0.3cm}\linkpatternunit=\unitsize

Here are a few examples in small size. Blank entries correspond to zeroes due to violation of the order.
\newcommand\vc[1]{\hbox{\vrule height 15pt depth 10pt width 0pt $\vcenter{\hbox{#1}}$}}
%\newcommand\vc[1]{\hbox{\vrule height 15pt depth 10pt width 0pt #1}}

%\begin{table}[h!]
%\caption{Table of $P^-$ for $N=4$, $K=2$ \rem{Brenti}}
%\centering
\begin{center}
Table of $P^-$ for $N=4$, $K=2$ \rem{Brenti}
\begin{tabular}{c|c|c|c|c|c|c|}
&\linkpattern[inverted]{1/1,4/4}&\linkpattern[inverted]{1/1,2/3,4/4}&\linkpattern[inverted]{1/1,3/4}&\linkpattern[inverted]{1/2,4/4}&\linkpattern[inverted]{1/2,3/4}&\linkpattern[inverted]{1/4,2/3}\\
&\path{-,-,+,+}&\path{-,+,-,+}&\path{-,+,+,-}&\path{+,-,-,+}&\path{+,-,+,-}&\path{+,+,-,-}\\
\hline
\vc{\path{-,-,+,+}}&$1$&$t^{-1}$&$0$&$0$&$0$&$t^{-2}$\\
\hline
\vc{\path{-,+,-,+}}&&$1$&$t^{-1}$&$t^{-1}$&$t^{-2}$&$t^{-1}$\\
\hline
\vc{\path{-,+,+,-}}&&&$1$&&$t^{-1}$&$0$\\
\hline
\vc{\path{+,-,-,+}}&&&&$1$&$t^{-1}$&$0$\\
\hline
\vc{\path{+,-,+,-}}&&&&&$1$&$t^{-1}$\\
\hline
\vc{\path{+,+,-,-}}&&&&&&$1$\\
\hline
\end{tabular}
\end{center}
%\end{table}

\smallskip

%\begin{table}[h!]
%\caption{Table of $P^+$ for $N=4$, $K=2$ \rem{LS}}
%\centering
\begin{center}
Table of $P^+$ for $N=4$, $K=2$ \rem{LS}
\begin{tabular}{c|c|c|c|c|c|c|}
&\linkpattern[inverted]{1/4,2/3}&\linkpattern[inverted]{1/2,3/4}&\linkpattern[inverted]{1/2,4/4}&\linkpattern[inverted]{1/1,3/4}&\linkpattern[inverted]{1/1,2/3,4/4}&\linkpattern[inverted]{1/1,4/4}\\
&\path{+,+,-,-}&\path{+,-,+,-}&\path{+,-,-,+}&\path{-,+,+,-}&\path{-,+,-,+}&\path{-,-,+,+}\\
\hline
\vc{\path{+,+,-,-}}&$1$&$t^{-1}$&$t^{-2}$&$t^{-2}$&$t^{-3}(1+t^2)$&$t^{-4}$\\
\hline
\vc{\path{+,-,+,-}}&&$1$&$t^{-1}$&$t^{-1}$&$t^{-2}$&$t^{-3}$\\
\hline
\vc{\path{+,-,-,+}}&&&$1$&&$t^{-1}$&$t^{-2}$\\
\hline
\vc{\path{-,+,+,-}}&&&&$1$&$t^{-1}$&$t^{-2}$\\
\hline
\vc{\path{-,+,-,+}}&&&&&$1$&$t^{-1}$\\
\hline
\vc{\path{-,-,+,+}}&&&&&&$1$\\
\hline
\end{tabular}
\end{center}
%\end{table}

The only non-monomial polynomial in $P^+$ corresponds to the
two Dyck strip decompositions
$\begin{tikzpicture}[x=\unitsize,y=\unitsize,baseline=0]
\tikzset{vertex/.style={circle,fill=black,inner sep=1.2pt,outer sep=0pt}}%
\tikzset{edge/.style={very thick}}%
\rawpath{+,+,-,-}\rawpath{-,+,-,+}\end{tikzpicture}$ 
and
$\begin{tikzpicture}[x=\unitsize,y=\unitsize,baseline=0]
\tikzset{vertex/.style={circle,fill=black,inner sep=1.2pt,outer sep=0pt}}%
\tikzset{edge/.style={very thick}}%
\rawpath{+,+,-,-}\rawpath{-,+,-,+}
\tikzset{vertex/.style={}}
\tikzset{edge/.style={}}
\rawpath{+,-,+,-}
\end{tikzpicture}$.

\rem{to go from one table to the other one, invert, transpose,
$t\to -t$ and reverse the order}

\section{An example of rule I}
\newcommand\exampletree[3]{%
\treestyle
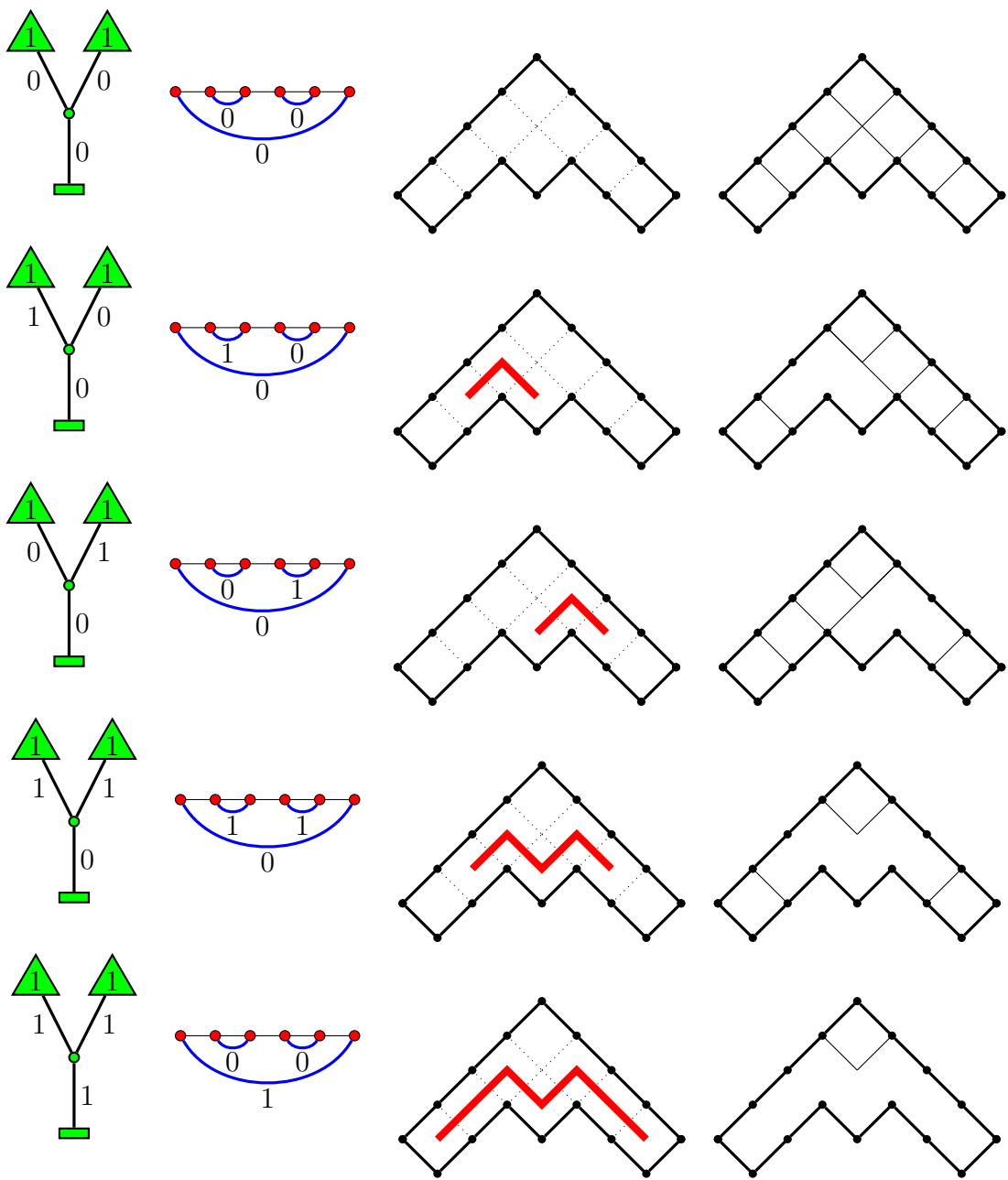
\begin{tikzpicture}[x=\unitsize,y=\unitsize]
\node[root] {} [tree]
%child 
%{
%  node[leaf] {0}
%  edge from parent node[left] {0}
%}
child 
{
  node[vertex] {}
  child 
  {
    node[leaf] {1}
    edge from parent node[left] {#1}
  }
  child 
  {
    node[leaf] {1}
    edge from parent node[right] {#2}
  }
  edge from parent node[right=-2pt] {#3}
}
%child 
%{
%  node[leaf] {0}
%  edge from parent node[right] {0}
%}
;
\end{tikzpicture}%
}
\newcommand\examplelinkpattern[3]{%
\raisebox{1.5cm}{%
\begin{tikzpicture}[x=\unitsize,y=\unitsize,baseline=0]
\tikzset{vertex/.style={circle,draw=black,fill=red,inner sep=1.5pt}}%
\tikzset{edge/.style={bend right=60,draw,very thick,blue}}%
\numberedfalse%
\linkpattern[tikzstarted]{1/6,2/3,4/5}%
\node at (2.5,-0.75) {#1};
\node at (4.5,-0.75) {#2};
\node at (3.5,-1.8) {#3};
\end{tikzpicture}}%
}

%\begin{figure}[h!]
\begin{center}
\setlength{\unitsize}{0.5cm}
\exampletree{0}{0}{0}
\quad
\examplelinkpattern{0}{0}{0}
\quad
\begin{tikzpicture}[x=\unitsize,y=\unitsize,baseline=0]
\tikzset{vertex/.style={circle,fill=black,inner sep=1.2pt,outer sep=0pt}}%
\tikzset{edge/.style={very thick}}%
\numberedfalse%
\rawpath{+,+,+,+,-,-,-,-}%
\rawpath{-,+,+,-,+,-,-,+}%
\tikzset{vertex/.style={}}%
\tikzset{edge/.style={dotted}}%
\rawpath{+,+,+,-,+,-,-,-}%
\rawpath{+,+,-,+,-,+,-,-}%
\rawpath{+,-,+,-,+,-,+,-}%
\end{tikzpicture}
\quad
\begin{tikzpicture}[x=\unitsize,y=\unitsize,baseline=0]
\tikzset{vertex/.style={circle,fill=black,inner sep=1.2pt,outer sep=0pt}}%
\tikzset{edge/.style={very thick}}%
\rawpath{+,+,+,+,-,-,-,-}%
\rawpath{-,+,+,-,+,-,-,+}%
\tikzset{vertex/.style={}}%
\tikzset{edge/.style={}}%
\rawpath{+,+,+,-,+,-,-,-}%
\rawpath{+,+,-,+,-,+,-,-}%
\rawpath{+,-,+,-,+,-,+,-}%
\end{tikzpicture}

\exampletree{1}{0}{0}
\quad
\examplelinkpattern{1}{0}{0}
\quad
\begin{tikzpicture}[x=\unitsize,y=\unitsize,baseline=0]
\tikzset{vertex/.style={circle,fill=black,inner sep=1.2pt,outer sep=0pt}}%
\tikzset{edge/.style={very thick}}%
\numberedfalse%
\rawpath{+,+,+,+,-,-,-,-}%
\rawpath{-,+,+,-,+,-,-,+}%
\tikzset{vertex/.style={}}%
\tikzset{edge/.style={dotted}}%
\rawpath{+,+,+,-,+,-,-,-}%
\rawpath{+,+,-,+,-,+,-,-}%
\rawpath{+,-,+,-,+,-,+,-}%
\draw[red, line width=3pt] (2,1) -- (3,2) -- (4,1);
\end{tikzpicture}
\quad
\begin{tikzpicture}[x=\unitsize,y=\unitsize,baseline=0]
\tikzset{vertex/.style={circle,fill=black,inner sep=1.2pt,outer sep=0pt}}%
\tikzset{edge/.style={very thick}}%
\numberedfalse%
\rawpath{+,+,+,+,-,-,-,-}%
\rawpath{-,+,+,-,+,-,-,+}%
\tikzset{vertex/.style={}}%
\tikzset{edge/.style={}}%
\rawpath{+,+,+,-,+,-,-,-}%
\rawpath{+,+,+,-,-,+,-,-}%
\rawpath{+,-,+,-,+,-,+,-}%
\end{tikzpicture}

\exampletree{0}{1}{0}
\quad
\examplelinkpattern{0}{1}{0}
\quad
\begin{tikzpicture}[x=\unitsize,y=\unitsize,baseline=0]
\tikzset{vertex/.style={circle,fill=black,inner sep=1.2pt,outer sep=0pt}}%
\tikzset{edge/.style={very thick}}%
\numberedfalse%
\rawpath{+,+,+,+,-,-,-,-}%
\rawpath{-,+,+,-,+,-,-,+}%
\tikzset{vertex/.style={}}%
\tikzset{edge/.style={dotted}}%
\rawpath{+,+,+,-,+,-,-,-}%
\rawpath{+,+,-,+,-,+,-,-}%
\rawpath{+,-,+,-,+,-,+,-}%
\draw[red, line width=3pt] (4,1) -- (5,2) -- (6,1);
\end{tikzpicture}
\quad
\begin{tikzpicture}[x=\unitsize,y=\unitsize,baseline=0]
\tikzset{vertex/.style={circle,fill=black,inner sep=1.2pt,outer sep=0pt}}%
\tikzset{edge/.style={very thick}}%
\numberedfalse%
\rawpath{+,+,+,+,-,-,-,-}%
\rawpath{-,+,+,-,+,-,-,+}%
\tikzset{vertex/.style={}}%
\tikzset{edge/.style={}}%
\rawpath{+,+,+,-,+,-,-,-}%
\rawpath{+,+,-,+,+,-,-,-}%
\rawpath{+,-,+,-,+,-,+,-}%
\end{tikzpicture}

\exampletree{1}{1}{0}
\quad
\examplelinkpattern{1}{1}{0}
\quad
\begin{tikzpicture}[x=\unitsize,y=\unitsize,baseline=0]
\tikzset{vertex/.style={circle,fill=black,inner sep=1.2pt,outer sep=0pt}}%
\tikzset{edge/.style={very thick}}%
\numberedfalse%
\rawpath{+,+,+,+,-,-,-,-}%
\rawpath{-,+,+,-,+,-,-,+}%
\tikzset{vertex/.style={}}%
\tikzset{edge/.style={dotted}}%
\rawpath{+,+,+,-,+,-,-,-}%
\rawpath{+,+,-,+,-,+,-,-}%
\rawpath{+,-,+,-,+,-,+,-}%
\draw[red, line width=3pt] (2,1) -- (3,2) -- (4,1) -- (5,2) -- (6,1);
\end{tikzpicture}%
\quad
\begin{tikzpicture}[x=\unitsize,y=\unitsize,baseline=0]
\tikzset{vertex/.style={circle,fill=black,inner sep=1.2pt,outer sep=0pt}}%
\tikzset{edge/.style={very thick}}%
\numberedfalse%
\rawpath{+,+,+,+,-,-,-,-}%
\rawpath{-,+,+,-,+,-,-,+}%
\tikzset{vertex/.style={}}%
\tikzset{edge/.style={}}%
\rawpath{+,+,+,-,+,-,-,-}%
\rawpath{+,-,+,-,+,-,+,-}%
\end{tikzpicture}%

\exampletree{1}{1}{1}
\quad
\examplelinkpattern{0}{0}{1}
\quad
\begin{tikzpicture}[x=\unitsize,y=\unitsize,baseline=0]
\tikzset{vertex/.style={circle,fill=black,inner sep=1.2pt,outer sep=0pt}}%
\tikzset{edge/.style={very thick}}%
\numberedfalse%
\rawpath{+,+,+,+,-,-,-,-}%
\rawpath{-,+,+,-,+,-,-,+}%
\tikzset{vertex/.style={}}%
\tikzset{edge/.style={dotted}}%
\rawpath{+,+,+,-,+,-,-,-}%
\rawpath{+,+,-,+,-,+,-,-}%
\rawpath{+,-,+,-,+,-,+,-}%
\draw[red, line width=3pt] (1,0) -- (3,2) -- (4,1) -- (5,2) -- (7,0);
\end{tikzpicture}%
\quad
\begin{tikzpicture}[x=\unitsize,y=\unitsize,baseline=0]
\tikzset{vertex/.style={circle,fill=black,inner sep=1.2pt,outer sep=0pt}}%
\tikzset{edge/.style={very thick}}%
\numberedfalse%
\rawpath{+,+,+,+,-,-,-,-}%
\rawpath{-,+,+,-,+,-,-,+}%
\tikzset{vertex/.style={}}%
\tikzset{edge/.style={}}%
\rawpath{+,+,+,-,+,-,-,-}%
\end{tikzpicture}%
\end{center}
%\caption{blah}
%\end{figure}

\section{A more explicit formula for $(P^+)^{-1}$}
Throughout this section, we again identify paths and binary strings.
We describe the set  $\mathcal{L}(\beta)$, which is the ``transposed'' 
set of  $\mathcal{F}(\beta)$.   

{\it A  linkage $w$} of a path $\beta$ is a set of pairs of integers 
from $[N]:=\{1,\ldots,N\}$ satisfying:
\begin{enumerate}
\item Each integer in $N$ appears exactly once in $w$. 
\item If a pair $(i,j)\in w, i,j\in [N]$ , then $\beta_i=1$ and $\beta_j=2$. 
\item Suppose $i$ and $j$, ($i<j$)  are paired. Then, there is no 
pair of $k$ and $l$ ($k<l$) such that $i<k<j<l$ or $k<i<l<j$. 
\end{enumerate}
Note that there are several linkages for a given path $v$, however, 
we recover a path from a given linkage.  
\begin{defn}
$\Cap(\beta)$ is a set of all possible linkages of the path $\beta$.  
\end{defn}

We need some terminology for pairs to define a map from an element 
$w\in\Cap(\beta)$ to the set of paths $\mathcal{L}(\beta)$.
\begin{enumerate}
\item A pair $(i,j)$ is said to an {\it ordered}  
(resp. {\it reversed}) pair if $i<j$ (resp. $i>j$).  
\item A pair of $k$ and $l$, $k<l$, i.e., a pair $(k,l)$ or $(l,k)$, is said 
to be {\it inside} of a pair of $i$ and $j$ if $i<k<l<j$, where 
$i,j,k,l\in[N]$.  
\end{enumerate}
%We call a reversed pair as a {\it r-pair}. 

We define an operation {\it r-flip} acting on a reversed pair $P$ 
in a linkage $w$ as follows. 
We flip $i$ and $j$ in $P$, all reversed pairs inside of $P$ and  keep all 
ordered pairs unchanged.  

\begin{defn}
$\mathcal{L}'(\beta;w)$ be the all possible paths recovered 
from linkages obtained by r-flipping (or without r-flipping) 
the linkage $w\in\Cap(\beta)$. 
\end{defn}

\begin{defn}
The set of paths by taking the union of $\mathcal{L}'(\beta;w)$ 
with respect to $w$: 
\begin{eqnarray*}
\mathcal{L}(\beta):=\bigcup_{w\in\Cap(\beta)}\mathcal{L}'(\beta;w).
\end{eqnarray*}
\end{defn}

In general, $\mathcal{L}'(\beta;w)\cap\mathcal{L}'(\beta;w')\neq\emptyset$ 
if $w,w'\in\Cap(\beta)$. 
Let $\alpha,\beta$ be two paths and $\alpha\in\mathcal{L}(\beta)$. 
The function $d(\alpha,\beta)$ depends only on the two paths, and this function 
counts the number of flipped r-pairs in $w\in\Cap(\beta)$ to obtain 
the path $\alpha$. 
Therefore, the number of flipped r-pairs to obtain $\alpha$ from $\beta$ are 
independent of the choice of a linkage. 

It is not hard to see
that the set $\mathcal{L}(\beta)$ describes exactly the set of $\alpha$
for which the summand in the formula of Corollary \ref{inverse-LS} is non-zero.
Therefore we have the slightly more explicit formula:
\[
m_\beta
 =
\sum_{\alpha\in\mathcal{L}(\beta)}(-t)^{-d(\alpha,\beta)}C^+_{\alpha} 
\]
%
% to fix url issues with line breaking
\let\oldurl\url
\makeatletter
\renewcommand*\url{%
        \begingroup
        \let\do\@makeother
        \dospecials
        \catcode`{1
        \catcode`}2
        \catcode`\ 10
        \url@aux
}
\newcommand*\url@aux[1]{%
        \setbox0\hbox{\oldurl{#1}}%
        \ifdim\wd0>\linewidth
                \strut
                \\
                \vbox{%
                        \hsize=\linewidth
                        \kern-\lineskip
                        \raggedright
                        \strut\oldurl{#1}%
                }%
        \else
                \hskip0pt plus\linewidth
                \penalty0
                \box0
        \fi
        \endgroup
}
\makeatother
% to fix MR issues
\gdef\MRshorten#1 #2MRend{#1}%
\gdef\MRfirsttwo#1#2{\if#1M%
MR\else MR#1#2\fi}
\def\MRfix#1{\MRshorten\MRfirsttwo#1 MRend}
\renewcommand\MR[1]{\relax\ifhmode\unskip\spacefactor3000 \space\fi
  \MRhref{\MRfix{#1}}{{\tiny\MRfix{#1}}}}
\renewcommand{\MRhref}[2]{%
 \href{http://www.ams.org/mathscinet-getitem?mr=#1}{#2}}
\DeclareUrlCommand\path{\urlstyle{tt}}%cause of conflict with my macros
\bibliographystyle{amsplainhyper} 
\bibliography{biblio}
\end{document}